\documentclass[10pt]{amsart}
\usepackage{amsmath,amsfonts,amssymb,amsthm,latexsym,amscd}
\usepackage[all]{xy}
\usepackage[mathscr]{eucal}
\theoremstyle{plain}
\newtheorem{thm}{Theorem}[section]
\newtheorem{prop}[thm]{Proposition}
\newtheorem{cor}[thm]{Corollary}
\newtheorem{lem}[thm]{Lemma}

\theoremstyle{definition}
\newtheorem{defn}[thm]{Definition}

\newtheorem{expl}[thm]{Example}

\theoremstyle{remark}
\newtheorem{rem}[thm]{Remark}

\numberwithin{equation}{section}

\newcommand{\LL}{\mathbf L}

\newcommand{\la}{\longrightarrow}

\newcommand{\sA}{\mathcal A}
\newcommand{\sB}{\mathcal B}
\newcommand{\sC}{\mathcal C}
\newcommand{\sD}{\mathcal D}

\newcommand{\sH}{\mathcal H}

\newcommand{\sK}{\mathcal K}
\newcommand{\sL}{\mathcal L}

\newcommand{\sP}{\mathcal P}
\newcommand{\sQ}{\mathcal Q}
\newcommand{\sR}{\mathcal R}

\newcommand{\sW}{\mathcal W}
\newcommand{\sO}{\mathcal O}

\newcommand{\bL}{\mathbf L}

\newcommand{\bH}{\mathbf H}

\newcommand{\RR}{\mathbb R}
\newcommand{\QQ}{\mathbb Q}
\newcommand{\ZZ}{\mathbb Z}

\newcommand{\id}{\textup{id}}

\newcommand{\lra}{\longrightarrow}
\newcommand{\co}{\colon\!}

\newcommand{\holim}{\textup{holim}}
\newcommand{\hocolim}{\textup{hocolim}}
\newcommand{\colim}{\textup{colim}}

\newcommand{\ko}{\textit{M}}

\newcommand{\pt}{\textup{pt}}

\newcommand{\TOP}{\textup{TOP}}
\newcommand{\STOP}{\textup{STOP}}
\newcommand{\Or}{\textup{O}}
\newcommand{\SOr}{\textup{SO}}
\newcommand{\PL}{\textup{PL}}
\newcommand{\SPL}{\textup{SPL}}

\newcommand{\smin}{\smallsetminus}

\newcommand{\ee}{\'{e} }
\newcommand{\ds}{\displaystyle}

\newcommand{\colimsub}[1]{\begin{array}[t]{cc} \textup{colim} \\
[-1.7mm] \scriptstyle{#1} \end{array}}

\newcommand{\hocolimsub}[1]{\begin{array}[t]{cc} \textup{hocolim} \\
[-1.7mm] \scriptstyle{#1} \end{array}}

\newcommand{\holimsub}[1]{\begin{array}[t]{cc} \textup{holim} \\
[-1.3mm] \scriptstyle{#1} \end{array}}

\newcommand{\twosub}[2]{\begin{array}{cc}
\scriptstyle{#1} \\  [-1mm] \scriptstyle{#2}  \end{array}}

\begin{document}

\title[Pontryagin classes]{On the construction and topological invariance of the Pontryagin classes}

\author{Andrew Ranicki and Michael Weiss}

\address{School of Mathematics \\ University of Edinburgh
\\ Edinburgh EH9 3JZ \\ Scotland, UK}
\email{a.ranicki@ed.ac.uk}

\address{Dept. of Math. Sciences \\ University of Aberdeen
\\ Aberdeen AB24 3UE \\ Scotland, UK}
\email{mweiss@maths.abdn.ac.uk}

\begin{abstract} We use sheaves and algebraic $L$-theory to construct the rational Pontryagin classes
of fiber bundles with fiber $\RR^n$. This amounts to an alternative proof of Novikov's
theorem on the topological invariance of the rational Pontryagin classes of vector bundles.
Transversality arguments and torus tricks are avoided. \end{abstract}

\date{January 7, 2009}

\maketitle

\section{Introduction}
The ``topological invariance of the rational Pontryagin classes'' was
originally the statement that for a \emph{homeomorphism} of smooth
manifolds, $f\co N\to N'$, the induced map
\[ f^*\co H^*(N';\QQ)\to H^*(N;\QQ) \]
takes the Pontryagin classes $p_i(TN')\in H^{4i}(N';\QQ)$ of the
tangent bundle of $N'$ to the Pontryagin classes $p_i(TN)$ of the
tangent bundle of $N$.  The topological invariance was proved by
Novikov \cite{Novikov}, some 40 years ago.  This breakthrough result
and its torus-related method of proof have stimulated many subsequent
developments in topological manifolds, notably the formulation of the
Novikov conjecture, the Kirby-Siebenmann structure theory and the
Chapman-Ferry-Quinn et al.  controlled topology.

\medskip
The topological invariance of the rational Pontrjagin classes
has been reproved several times (Gromov
\cite{Gromov}, Ranicki \cite{RanickiNov}, Ranicki and Yamasaki
\cite{RanickiYamasaki}).  In this paper we give yet another proof,
using sheaf-theoretic ideas. We associate to a topological manifold $M$
a kind of ``tautological'' co-sheaf on $M$ of symmetric Poincar\'e chain complexes, in such a way that the
cobordism class is a topological invariant by construction. We then produce
excision and homotopy invariance theorems for the cobordism groups of such
co-sheaves, and use the Hirzebruch signature theorem
to extract the rational Pontryagin classes of $M$ from the cobordism
class of the tautological co-sheaf on $M$.

\begin{defn} Write $\TOP(n)$ for the space of homeomorphisms from $\RR^n$ to $\RR^n$
and $\PL(n)$ for the space (geometric realization of a simplicial set) of
PL-homeomorphisms from $\RR^n$ to $\RR^n$. Let $\TOP=\bigcup_n
\TOP(n)$ and $\PL=\bigcup_n \PL(n)$.
\end{defn}

The topological invariance can also be formulated in terms of
the classifying spaces, as the statement that the
homomorphism $H^*(B\TOP;\QQ)\to H^*(B\Or;\QQ)$ induced by the
inclusion $\Or\to \TOP$ is onto.\footnote{The point is that
topological $n$-manifolds have topological tangent bundles with
structure group $\TOP(n)$ which are sufficiently natural under
homeomorphisms. See \cite{Kister}.} Here
\[H^*(B\Or;\QQ)=\QQ[p_1,p_2,p_3,\dots] \] where $p_i\in
H^{4i}(B\Or;\QQ)$ is the $i$-th ``universal'' rational Pontryagin class.
We note also that the inclusion $B\SOr\to B\Or$ induces an isomorphism in
rational cohomology, and the inclusion $B\STOP\to B\TOP$ induces a surjection
in rational cohomology by a simple transfer argument. Therefore it is enough to
establish surjectivity of $H^*(B\STOP;\QQ)\to H^*(B\SOr;\QQ)$.

\medskip
Hirzebruch's signature theorem \cite{Hirz} expressed the signature
of a closed smooth oriented $4i$-dimensional manifold $M$ as the
evaluation on the fundamental class $[M] \in H_{4i}(M)$ of the
\emph{$\sL$-genus}
\[  \sL(TM)\in H^{4*}(M;\QQ), \]
that is
\[ \text{signature}(M)~=~ \langle \sL(TM),[M] \rangle \in \ZZ \subset \QQ~.\]
It follows that, for any closed smooth oriented $n$-dimensional manifold $N$
and a closed $4i$-dimensional framed submanifold $M \subset N \times
\RR^k$ (where \emph{framed} refers to a trivialized normal
bundle),
\[
{\rm signature}(M)~=~\langle \sL(TN),[M] \rangle \in \ZZ \subset
\QQ~. \] By Serre's finiteness theorem for homotopy groups,
\[ H^{n-4i}(N;\QQ)~\cong~\varinjlim_k~[\Sigma^kN_+,S^{n-4i+k}]\otimes \QQ \]
with $N_+=N \cup \{ \text{pt.}\}$, and by Pontryagin-Thom theory
$[\Sigma^kN_+,S^{n-4i+k}]$ can be identified with the bordism group
of closed framed $4i$-dimensional submanifolds $M \subset N \times
\RR^k$. It is thus possible to identify the component of $\sL(TN)$
in
\[ H^{4i}(N;\QQ)~\cong~\hom(H^{n-4i}(N;\QQ),\QQ)\]
with the linear map
\[ H^{n-4i}(N;\QQ)\lra \QQ~;~M \mapsto \text{signature}(M)~.\]
The assumption that $N$ be closed can be discarded if we use
cohomology with compact supports where appropriate. Then we identify
the component of $\sL(TN)$ in
\[ H^{4i}(N;\QQ)~\cong~\hom(H^{n-4i}_c(N;\QQ),\QQ)\]
with the linear map
\[ H^{n-4i}_c(N;\QQ)\to \QQ~;~M \mapsto \text{signature}(M)~.\]
Now we can choose $N$ in such a way that the classifying map
$\kappa:N\to B\SOr(n)$ for the tangent bundle is highly connected,
say $(4i+1)$-connected. Then $\kappa^*:H^{4i}(B\SOr(n);\QQ)\to
H^{4i}(N;\QQ)$ takes $\sL$ of the universal oriented $n$-dimensional
vector bundle to $\sL(TN)$, and so the above description of
$\sL(TN)$ in terms of signatures can be taken as a \emph{definition}
of the universal $\sL\in H^{4i}(B\SOr(n);\QQ)$, or even $\sL\in
H^{4i}(B\SOr;\QQ)$.

\smallskip
Since this definition relies almost exclusively on transversality
arguments, which carry over to the PL setting, we can deduce
immediately, as Thom did, that Hirzebruch's $\sL$-genus extends to
\[ \sL\in H^{4*}(B\SPL;\QQ)~. \]
As the rational Pontryagin classes are polynomials (with rational
coefficients) in the components of $\sL$, the PL invariance of the
rational Pontryagin classes follows from this rather straightforward
argument. It was also clear that the \emph{topological} invariance
of the rational Pontryagin classes would follow from an appropriate
transversality theorem in the setting of topological manifolds.
However, Novikov's proof did not exactly deduce the topological
invariance of the Pontryagin classes from a topological
transversality statement. Instead, he proved that signatures of the
submanifolds were homeomorphism invariants by showing that for a
homeomorphism $f:N \to N'$ of closed \emph{smooth} oriented
$n$-dimensional manifolds $N,N'$ and a closed framed
$4i$-dimensional submanifold $M' \subset N' \times \RR^k$ it is
possible to make the proper map $f \times \id : N \times \RR^k \to
N' \times \RR^k$ transverse regular at $M'$, keeping it proper, and
the smooth transverse image $M \subset N \times \RR^k$ has
\[\text{signature}(M)~=~\text{signature}(M') \in \ZZ~. \]
This was done using non-simply-connected methods.
Subsequently, it was found that the ideas in Novikov's proof could be extended
and combined with non-simply-connected surgery theory to prove
transversality for topological (non-smooth, non-PL) manifolds.  Details
on that can be found in \cite{KirbySiebenmann}.  It is now known that
$H^*(B\TOP;\QQ)\cong H^*(B\Or;\QQ)$.
\medskip

The collection of the signatures of framed $4i$-dimensional
submanifolds $M \subset N \times \RR^k$ of a topological
$n$-dimensional manifold $N$ was generalized in Ranicki
\cite{RanickiTop} to a fundamental class $[N]_{\LL} \in
H_n(N;\LL^{\bullet})$ with coefficients in a spectrum
$\LL^{\bullet}$ of symmetric forms over $\ZZ$.  However,
\cite{RanickiTop} made some use of topological transversality. This
paper will remedy this by expressing $[N]_{\LL}$ as the cobordism
class of the tautological co-sheaf on $N$, using the local
Poincar\'e duality properties of $N$.

\medskip
In any case the ``modern'' point of view in the matter of the
topological invariance of rational Pontryagin classes is that it
merits a treatment separate from transversality discussions.
Topological manifolds ought to have (tangential) rational Pontryagin
classes because they satisfy a local form of Poincar\ee
duality. More precisely, if $N$ is a topological $n$-manifold, then
for any open set $U\subset N$ we have a Poincar\'e duality isomorphism
between the homology of $U$ and the cohomology of $U$ with compact
supports. The task is then to use this refined form of Poincar\ee
duality to make invariants in the homology or cohomology of $N$.
Whether or not these invariants can be expressed as characteristic
classes of the topological tangent bundle of $N$ becomes a question
of minor importance.

\bigskip
\emph{Remark.} In this paper we make heavy use of homotopy direct and homotopy inverse limits of diagrams
of chain complexes and chain maps. They can be defined like homotopy direct and homotopy inverse limits of
diagrams of spaces. That is to say, the Bousfield-Kan formulae \cite{BousfieldKan} for homotopy direct and homotopy
inverse limits of diagrams of spaces can easily be adapted to diagrams of chain complexes:
products (of spaces) should be replaced by tensor products (of chain complexes), and where standard
simplices appear they should, as a rule, be replaced by their cellular chain complexes. \newline
We often rely on \cite[\S9]{DwyerKan}, a collection of conversion and comparison theorems
for homotopy direct and homotopy inverse limits.

\section{Duality and $L$-theory: Generalities}
\label{sec-easysetting} In the easiest setting, we start with an
additive category $\sA$ and the category $\sC$ of all chain
complexes in it, graded over $\ZZ$ and bounded from above and below.
We assume given a functor
\[  (C,D)\mapsto C\boxtimes D \]
from $\sC\times\sC$ to chain complexes of abelian groups. This is
subject to bilinearity and symmetry conditions:
\begin{itemize}
\item for fixed $D$ in $\sC$, the functor $C\mapsto C\boxtimes D$ takes contractible objects
to contractible objects, and takes homotopy cocartesian (= homotopy
pushout) squares to homotopy cocartesian squares~;
\item there is a binatural isomorphism $\tau\co C\boxtimes D\to D\boxtimes C$
satisfying $\tau^2=\id$.
\end{itemize}
We assume that every object $C$ in $\sC$ has a ``dual''
(w.r.t.~$\boxtimes$). This means that the functor
\[   D\mapsto H_0(C\boxtimes D) \]
is co-representable in the homotopy category $\sH\sC$, so that there
exists $C^{-*}$ in $\sC$ and a natural isomorphism
\[  H_0(C\boxtimes D)\cong[C^{-*},D] \]
where the square brackets denote chain homotopy classes of maps.
Let's note that in this case the $n$-fold suspension
$\Sigma^nC^{-*}$ co-represents the functor $C\mapsto H_n(C\boxtimes
D)$ in $\sH\sC$.

\begin{expl} $\sA$ is the category of f.g. free left modules over $\ZZ[\pi]$,
for a fixed group $\pi$. For $C$ and $D$ in $\sC$ we let
\[  C\boxtimes D= C^t\otimes_{\ZZ[\pi]}D~, \]
using the standard involution $\sum a_g g\mapsto \sum a_g g^{-1}$ on
$\ZZ[\pi]$ to turn $C$ into a chain complex of right
$\ZZ[\pi]$-modules $C^t$. Then the dual $C^{-*}$ of any $C$ in $\sC$
exists and can be defined explicitly as the chain complex of right
module homomorphisms $\hom_{\ZZ[\pi]}(C^t,\ZZ[\pi])$ on which
$\ZZ\pi$ acts by left multiplication:
\[  (rf)(c):=r(f(c))  \]
for $r\in \ZZ[\pi]$, $c\in C^t$ and
$f\in\hom_{\ZZ[\pi]}(C^t,\ZZ[\pi])$.
\end{expl}

\begin{defn} An $n$-dimensional
symmetric algebraic Poincar\ee complex in $\sC$ consists of an
object $C$ in $\sC$ and an $n$-dimensional cycle $\varphi$ in
$(C\boxtimes C)^{h\ZZ/2}$ whose image in $H_n(C\boxtimes C)$ is
nondegenerate (i.e., adjoint to a homotopy equivalence
$\Sigma^nC^{-*}\to C$).
\end{defn}

\begin{rem} $(C\boxtimes
C)^{h\ZZ/2}=\hom_{\ZZ[\ZZ/2]}(W,C\boxtimes C)$ where $W$ is your
favorite projective resolution of the trivial module $\ZZ$ over the
ring $\ZZ[\ZZ/2]$.
\end{rem}

\medskip
\emph{Example}: Suppose that $\sA$ is the category of f.g. free
abelian groups. Let $C$ be the cellular chain complex of a space $X$
which is the realization of a f.g. simplicial set, and also a
Poincar\ee duality space of formal dimension $n$. Then you can use
an Eilenberg-Zilber diagonal
\[  W\otimes C \lra C\otimes C \]
(respecting $\ZZ/2$-actions) and evaluate the adjoint $C\to
(C\otimes C)^{h\ZZ/2}$ on a fundamental cycle to get a nondegenerate
$n$-cycle $\varphi\in (C\otimes C)^{h\ZZ/2}$.

\begin{defn} An $(n+1)$-dimensional symmetric algebraic Poincar\ee pair in $\sC$ consists of
a morphism $f\co C\to D$ in $\sC$, an $n$-dimensional cycle
$\varphi$ in $(C\boxtimes C)^{h\ZZ/2}$ and an $(n+1)$-dimensional
chain $\psi$ in $(D\boxtimes D)^{h\ZZ/2}$ such that $(f\boxtimes
f)(\varphi)=
\partial\psi$ and
\begin{itemize}
\item the image of $\varphi$ in $H_n(C\boxtimes C)$ is nondegenerate,
\item the image of $\psi$ in $H_{n+1}(D\boxtimes D/C)$ is nondegenerate, where
$D/C$ is shorthand for the algebraic mapping cone of $f\co C\to D$.
\end{itemize}
The \emph{boundary} of the SAP pair $(C\to D,\psi,\varphi)$ is
$(C,\varphi)$, an $n$-dimensional SAPC.
\end{defn}

With these definitions, it is straightforward to design bordism
groups $L^n(\sC)$ of $n$-dimensional SAPCs in $\sC$. Elements of
$L^n(\sC)$ are represented by $n$-dimensional SAPCs in $\sC$. We say that
two such representatives, $(C,\varphi)$ and $(C',\varphi')$, are
\emph{bordant} if $(C\oplus C',\varphi\oplus-\varphi')$ is the
boundary of an SAP pair of formal dimension $n+1$. \newline
Furthermore there exist generalizations of the definitions of SAPC
and SAP pair (going in the direction of $m$-ads) which lead
automatically to the construction of a spectrum $\bL^\bullet(\sC)$
such that $\pi_n\bL^\bullet(\sC)\cong L^n(\sC)$. These activities
come under the heading \emph{symmetric $L$-theory of $\sC$}. This is
an idea going back to Mishchenko \cite{Mishchenko}. It was pointed
out in \cite{RanickiLms} that there is an analogue of Mishchenko's
setup, the \emph{quadratic $L$-theory of $\sC$}, where ``homotopy
fixed points'' of $\ZZ/2$-actions are replaced by ``homotopy
orbits''. For example:

\begin{defn} An $n$-dimensional
quadratic algebraic Poincar\ee complex in $\sC$ consists of an
object $C$ in $\sC$ and an $n$-dimensional cycle $\varphi$ in
$(C\boxtimes C)_{h\ZZ/2}$ whose image in $H_n(C\boxtimes C)$ (under
the transfer) is nondegenerate.
\end{defn}

\begin{rem} $(C\boxtimes
C)_{h\ZZ/2}=W\otimes_{\ZZ[\ZZ/2]}(C\boxtimes C)$ where $W$ is that
resolution, as above.
\end{rem}

\medskip
The bordism groups of QAPCs are denoted $L_n(\sC)$ and the
corresponding spectrum is $\bL_\bullet(\sC)$. There is a
(norm-induced) comparison map $\bL_\bullet(\sC)
\to\bL^\bullet(\sC)$. On the algebraic side, a key difference
between $L_n(\sC)$ and $L^n(\sC)$ is that elements of $L_n(\sC)$ can
always be represented by ``short'' chain complexes (concentrated in
degree $k$ if $n=2k$, and in degrees $k$ and $k+1$ if $n=2k+1$),
while that is typically not the case for elements of $L^n(\sC)$. On
the geometric side, $L_n(\sC)$ also has the immense advantage of
being directly relevant to differential topology as a surgery
obstruction group, for the right choice(s) of $\sC$. But the
comparison maps $L_n(\sC)\to L^n(\sC)$ are always isomorphisms away
from the prime 2, and since we are interested in rational questions,
there is no need to make a very careful distinction between
symmetric and quadratic $L$-theory here.

\begin{expl} If $\sA$ is the category of f.g. free abelian groups, $\sC$ the corresponding
chain complex category, then
\[ L_n(\sC)\cong\left\{\begin{array}{lcl}
\ZZ && n\equiv 0 \mod 4 \\
0  && n\equiv 1 \mod 4 \\
\ZZ/2 && n\equiv 2 \mod 4 \\
0 && n\equiv 3 \mod 4
\end{array} \right.
\]
If $\sA$ is the category of f.d. vector spaces over $\QQ$, and $\sC$
the corresponding chain complex category, then
\[ L_n(\sC)\cong L^n(\sC) \cong\left\{\begin{array}{lcl}
\ZZ\oplus(\ZZ/2)^\infty\oplus(\ZZ/4)^\infty && n\equiv 0 \mod 4 \\
0  && n\equiv 1,2,3 \mod 4
\end{array} \right.
\]
where $(...)^\infty$ denotes a countably infinite direct sum.
\end{expl}

\smallskip
\begin{rem}
\label{rem-triangulated}
We need a mild generalization of the setup above. Again we start
with an additive category $\sA$. Write $\sB(\sA)$ for the category of all
chain complexes of $\sA$-objects, bounded from below (but not necessarily from above). We suppose
that a full subcategory $\sK$ of $\sB(\sA)$ has been specified, closed
under suspension, desuspension, homotopy equivalences, direct sums
and mapping cone constructions, so that the homotopy category $\sH\sK$ is a triangulated
subcategory of $\sH\sB(\sA)$. We assume given a functor
\[  (C,D)\mapsto C\boxtimes D \]
from $\sK\times\sK$ to chain complexes of abelian groups. This is
subject to the usual bilinearity and symmetry conditions:
\begin{itemize}
\item for fixed $D$ in $\sK$, the functor $C\mapsto C\boxtimes D$ takes contractible
objects to contractible ones
and preserves homotopy cocartesian (= homotopy
pushout) squares;
\item there is a binatural isomorphism $\tau\co C\boxtimes D\to D\boxtimes C$
satisfying $\tau^2=\id$.
\end{itemize}
We assume that every object $C$ in $\sK$ has a ``dual''
(w.r.t.~$\boxtimes$). This means that the functor $D\mapsto
H_0(C\boxtimes D)$ on $\sH\sK$ is co-representable. From these data
we construct $L$-theory spectra $\bL_\bullet(\sK)$ and
$\bL^\bullet(\sK)$ as before. (Some forward ``hints'': Our choice
of additive category $\sA$ is fixed from definition~\ref{defn-A} onwards,
and for $\sK$ we take the category $\sD'$ defined in section~\ref{sec-zoo}.)
\end{rem}

\section{Chain complexes in a local setting}
\label{sec-introlocal}

Let $X$ be a locally compact, Hausdorff and separable space. Let $\sO(X)$ be
the poset of open subsets of $X$. We introduce an additive category
$\sA=\sA_X$ whose objects are free abelian groups (typically not
finitely generated) equipped with a system of subgroups indexed by
$\sO(X)$.

\begin{defn}
\label{defn-A}
An object of $\sA$ is a free abelian group $F$ with a basis $S$, together with
subgroups $F(U)\subset F$, for $U\in \sO(X)$, such that the following
conditions are satisfied.
\begin{itemize}
\item $F(\emptyset)=0$ and $F(X)=F$.
\item Each subgroup $F(U)$ is generated by a subset of $S$.
\item For $U,V\in \sO(X)$ we have $F(U\cap V)=F(U)\cap F(V)$.
\end{itemize}
A morphism in $\sA$ from $F_0$ to $F_1$ is a group homomorphism
$F_0\to F_1$ which, for every $U\in\sO(X)$, takes $F_0(U)$ to $F_1(U)$.
\end{defn}

\begin{expl}
\label{expl-singulargroup}
For $i\ge 0$, the $i$-th chain group $C_i=C_i(X)$ of the singular chain complex of $X$
has a preferred structure of an object of $\sA$. The preferred basis is the set $S_i=S_i(X)$ of singular
$i$-simplices in $X$. For open $U$ in $X$, let $C_i(U)\subset C_i(X)$ be the subgroup
generated by the singular simplices with image contained in $U$. The boundary operator
from $C_i(X)$ to $C_{i-1}(X)$ is an example of a morphism in $\sA$.
\end{expl}

\begin{defn} We write $\sB(\sA)$
for the category of chain complexes in $\sA$, graded over $\ZZ$
and bounded from below.
\end{defn}

\medskip
Next we list some conditions which we might impose on objects in
$\sB(\sA)$, to define a subcategory in which we can successfully do
$L$-theory. The kind of object that we are most interested in is
described in the following example.

\begin{expl}
\label{expl-guide} Take a map $f\co Y\to X$, where $Y$ is a
compact ENR (euclidean neighborhood retract). Let $C(f)$ be the
object of $\sB(\sA)$ defined as follows: $C(f)(X)$ is the singular chain
complex of $Y$, with the standard graded basis, and $C(f)(U)\subset C(f)(X)$ for $U\in \sO(X)$ is the
subcomplex generated by the singular simplices of $Y$ whose image is in $f^{-1}(U)$.
\end{expl}

We start by listing some of the obvious but
remarkable properties which we see in this example~\ref{expl-guide}.
Let $C$ be any object of $\sB(\sA)$.

\begin{defn}
\label{defn-sheaflist}
We say that $C$ satisfies the \emph{sheaf type
condition} if, for any subset $\sW$ of $\sO(X)$, the inclusion
\[  \sum_{V\in \sW} C(V) \la C(\bigcup_{V\in\sW}V) \]
is a homotopy equivalence.\footnote{The sum sign is for an internal sum taken in the chain complex $C(X)$, not an abstract
direct sum.}\newline
We say that $C$ satisfies \emph{finiteness condition} (i) if the following holds.
There exists an integer $a\ge 0$ such that
\begin{itemize}
\item[{}] for open sets $V_1\subset V_2$ in $X$ such that
the closure of $V_1$ is contained in $V_2$, the inclusion
$C(V_1)\to C(V_2)$ factors up to chain homotopy through a
chain complex $D$ of f.g. free abelian groups, with $D_i=0$ if $|i|>a$.
\end{itemize}
We say that $C$ satisfies \emph{finiteness condition} (ii) if
\begin{itemize}
\item[{}] there exists a compact $K$ in $X$ such that $C(U)$
depends only on $U\cap K$. (Then we say that $C$ is \emph{supported in
$K$}.)
\end{itemize}
\end{defn}

\begin{lem}
\label{lem-guide}
Example~\ref{expl-guide} satisfies the sheaf type
condition and the two finiteness conditions.
\end{lem}

\proof The sheaf type condition is well known from the standard
proofs of excision in singular homology. A crystal clear reference for this is
\cite[III.7.3]{Dold}. For finiteness condition
(i), choose a finite simplicial complex $Z$ and a retraction $r\co Z\to Y$
(with right inverse $j\co Y\to K$). Replacing the triangulation of $Z$
by a finer one if necessary, one can find a finite simplicial subcomplex $Z'$
of $Z$ containing $j(f^{-1}U)$ and such that $r(Z')\subset f^{-1}(V)$. Then the
inclusion $f^{-1}U\to f^{-1}(V)$ factors through $Z'$. The singular chain complex
of $Z'$ is homotopy equivalent to the cellular chain complex of $Z'$~, a
chain complex of f.g. free abelian groups which is zero in degrees $<0$
and in degrees $>\dim(Z)$.
Finiteness condition (ii) is satisfied with $K=f(Y)$. \qed

\begin{defn} Let $\sC\subset \sB(\sA)$ consist of the objects which satisfy the
sheaf type conditions and the two finiteness conditions.
\end{defn}

In the following proposition, we write $\sC_X$ and $\sC_Y$ etc. rather than $\sC$, to emphasize the dependence on
a space such as $X$ or $Y$.

\begin{defn} A map $f\co X\to Y$ induces a ``pushforward'' functor $f_*\co \sC_X\to \sC_Y$,
defined by $f_*C(U)=C(f^{-1}(U))$ for $C$ in $\sC_X$ and open $U\subset Y$.
\end{defn}

Returning to the shorter notation ($\sC$ for $\sC_X$), we spell out two elementary consequences
of the sheaf condition.

\begin{lem}
\label{lem-descent}
Let $C$ be an object of $\sC$ and let $\sW$ be a finite
subset of $\sO(X)$. If $\sW$ is closed under unions (i.e., for any $V,W\in \sW$ the
union $V\cup W$ is in $\sW$) then the inclusion-induced map
\[  C(\bigcap_{V\in\sW}V) \la \holimsub{V\in\sW} C(V) \]
is a homotopy equivalence.
\end{lem}

\proof
Choose $V_1,\dots,V_k\in \sW$ such that every element of $\sW$ is a union
of some (at least one) of the $V_i$~. For nonempty $S\subset \{1,\dots k\}$ let
$V_S=\bigcup_{i\in S}V_i$. There is an inclusion
\[ \holimsub{V\in\sW} C(V) \la \holimsub{\textup{nonempty }S\subset\{1,\dots,k\}} C(V_S)
\]
and we show first that this is a homotopy equivalence. It is induced
by a map $f$ of posets. On the right-hand side, we have the poset of nonempty subsets
of $\{1,\dots,k\}$ partially ordered by reverse inclusion, and on the left-hand side
we have $\sW$ itself, partially ordered by inclusion. The map $f$ is given
by $S\mapsto V_S$. Under these circumstances it is enough to show that $f$ has
a (right) adjoint $g$. But this is clear: for $V\in\sW$ let $g(V)=\{i\mid V_i\subset V\}$. \newline
Now it remains to show that the inclusion-induced map
\[
C(\,\bigcap_{i=1}^k V_i) \la \holimsub{\textup{nonempty }S\subset\{1,\dots,k\}} C(V_S)
\]
is a homotopy equivalence. We show this without any restrictive assumptions
on $V_1,\dots,V_k$. The map fits into a commutative square
\[
\xymatrix{
{C(\,\bigcap_{i=1}^{k-1} V_i\cap V_k)=C(\,\bigcap_{i=1}^k V_i)} \ar[r] \ar[d] &
{\holimsub{\emptyset\ne S\subset\{1,\dots,k\}} C(V_S)} \ar[d] \\
{\rule{0mm}{6mm}\holimsub{\emptyset \ne S\subset\{1,\dots,k-1\}} C(V_S\cap V_k)} \ar[r] &
{\rule{0mm}{6mm}\holimsub{\emptyset\ne S\subset\{1,\dots,k-1\}} \holim\left(C(V_S)\to C(V_{S\cup k})\leftarrow C(V_k)\right).}
}
\]
The left-hand vertical arrow in the square is a homotopy equivalence by inductive assumption.
The lower horizontal arrow is induced by homotopy equivalences
\[ C(V_S\cap V_k) \la \holim\left(C(V_S)\to C(V_{S\cup k})\leftarrow C(V_k)\right) \]
and is therefore also a homotopy equivalence. The right-hand vertical arrow is contravariantly
induced by a map of posets,
\[   (S,T) \mapsto T \]
where $S$ is a nonempty subset of $\{1,\dots,k-1\}$ and $T=S$ or $T=S\cup k$ or $T=k$.
Hence it is enough, by \cite[9.7]{DwyerKan}, to verify that the appropriate categorical ``fibers'' of this
map of posets have contractible classifying spaces.
For fixed $T'$, a nonempty subset of $\{1,\dots,k\}$, the appropriate fiber is the poset of all
$(S,T)$ as above with $T\subset T'$. It is easy to verify that the classifying space is
contractible. \qed

\begin{lem}
\label{lem-codescent}
Let $C$ be an object of $\sC$ and let $\sW$ be a
subset of $\sO(X)$. If $\sW$ is closed under intersections (i.e., for any $V,W\in \sW$ the
intersection $V\cap W$ is in $\sW$) then the inclusion-induced map
\[   \hocolimsub{V\in\sW} C(V) \la C(\bigcup_{V\in\sW}V)\]
is a homotopy equivalence.
\end{lem}

\proof In the case where $\sW$ is finite, the proof is analogous to that of lemma~\ref{lem-descent}.
The general case follows from the case where $\sW$ is finite by an obvious direct limit argument.
\qed

\section{A zoo of subcategories}
\label{sec-zoo}
The category $\sC=\sC_X$ defined in section~\ref{sec-introlocal} should be regarded as a provisional work environment.
It has two shortcomings.
\begin{itemize}
\item A morphism $f:C\to D$ in $\sC$ which induces homotopy equivalences $C(U)\to D(U)$ for
every open $U\subset X$ need not be a chain homotopy equivalence in $\sC$.
\end{itemize}
We can fix that rather easily, and will do so in this
section, by defining \emph{free} objects in $\sC$ and showing that all objects in $\sC$ have free resolutions.
This leads to a decomposition of $\sC$ into full subcategories $\sC'$ and $\sC''$, where $\sC'$ contains the
free objects and $\sC''$ contains the objects which we regards as \emph{weakly equivalent to } $0$.
\begin{itemize}
\item Given the decomposition of $\sC$ into $\sC'$ and $\sC''$, we are able to set up a good duality
theory either in $\sC'$ or, less formally, in $\sC$ modulo $\sC''$. The resulting quadratic
$L$-theory spectrum is still a functor of $X$, because $\sA,\sC,\sC',\sC''$ depend on $X$. For this functor we
are able to prove homotopy invariance, but not excision.
\end{itemize}
We solve this problem not by adding further conditions to the list in definition~\ref{defn-sheaflist},
but instead by defining a full subcategory $\sD$ of $\sC$ in terms of generators. To be more precise, $\sD$
is generated by all objects which are weakly equivalent to $0$ and all the examples of~\ref{expl-guide} obtained
from singular simplices $f\co \Delta^k\to X$, using the processes of extension, suspension and desuspension.
The good duality theory in $\sC'$ or $\sC$ modulo $\sC''$ restricts to a good duality theory in
$\sD'$ or $\sD$ modulo $\sD''$, where $\sD'=\sD\cap\sC'$ and $\sD''=\sD\cap \sC''$. The corresponding
$L$-theory functor $X\mapsto \bL^\bullet(\sD_X)$ satisfies homotopy invariance and excision. \newline
Unfortunately it is not clear that \emph{all} the objects of $\sC$ obtained as in example~\ref{expl-guide}
belong to $\sD$. They do however belong to $r\sD$, the idempotent completion of $\sD$ within $\sC$.
We are not able to prove excision for the functor $X\mapsto \bL^\bullet(r\sD_X)$, but we do have a long
exact ``Rothenberg'' sequence showing that the inclusion $\bL^\bullet(\sD_X)\to \bL^\bullet(r\sD_X)$ is a homotopy
equivalence away from the prime $2$.
\[
\xymatrix@M=10pt@W=10pt{
 & \sC' \ar@{>->}[d]  & \ar@{>->}[l] r\sD'  \ar@{>->}[d] & \ar@{>->}[l] \sD' \ar@{>->}[d] \\
\sB(A) & \ar@{>->}[l] \sC &  \ar@{>->}[l] r\sD &  \ar@{>->}[l] \sD \\
& \sC'' \ar@{>->}[u]  & \ar@{>->}[l] r\sD'' \ar@{>->}[u] & \ar@{>->}[l] \sD'' \ar@{>->}[u]
}
\]

\medskip
\begin{defn}
\label{defn-freeobject}
A nonempty open subset $U$ of $X$ determines a functor $F$ on $\sO(X)$ by
\[  F(V) = \left\{\begin{array}{cl} \ZZ & \textrm{ if } U\subset V \\ 0&
\textrm{ otherwise. }
\end{array} \right.
\]
With the obvious basis for $F(X)$, this becomes an object of $\sA$ which we call
\emph{free on one generator attached to $U$}. Any direct sum of such objects is called
\emph{free}.
\end{defn}

\begin{expl} The singular chain group $C_i(X)$ of $X$, with additional structure as
in example~\ref{expl-singulargroup}, is typically \emph{not} free in the sense of
\ref{defn-freeobject}.
\end{expl}

\begin{defn} Let $\sC'\subset \sC$ be the full subcategory consisting of the objects
which are free in every dimension. Let $\sC''\subset \sC$ be the full subcategory consisting of
the objects $C$ for which $C(U)$ is contractible, for all $U\in\sO(X)$.
\end{defn}

\begin{defn} A morphism $f\co C\to D$ in $\sC$ is a \emph{weak equivalence} if its mapping cone belongs to
$\sC''$.
\end{defn}

\begin{lem}
\label{lem-ortho}
Every morphism in $\sC$ from an object of $\sC'$ to an object
of $\sC''$ is nullhomotopic.
\end{lem}

\proof The nullhomotopy can be constructed by induction over skeleta,
using the following ``projective'' property of free objects in $\sA$.
Let a diagram
\[
\CD
 @.  B_0 \\
 @.  @VV f V \\
 A @>>> B_1
 \endCD
 \]
in $\sA$ be given where $f$ is strongly onto (that is, the induced map $B_0(U)\to B_1(U)$ is
onto for every $U$) and $A$ is free. Then there exists
$g\co A\to B_0$ making the diagram commutative. \qed

\begin{lem}
\label{lem-resolve}
For every object $D$ of $\,\sC$, there exists an object $C$ of
$\,\sC'$ and a morphism $C\to D$ which is a weak equivalence.
\end{lem}

\proof The morphism $C\to D$ can be constructed inductively using
the fact that, for every $B$ in $\sA$, there exists a free $A$ in $\sA$
and a morphism $A\to B$ which is strongly onto. \qed

\medskip
The next lemma means that objects of $\sC'$ are ``cofibrant'':

\begin{lem}
\label{lem-cofibrant} Let $f\co C\to D$ and $g\co E\to D$ be morphisms
in $\sC$. Suppose that $C$ is in $\sC'$ and $g$ is a weak equivalence.
Then there exists a morphism $f^\sharp\co C\to E$ such that
$gf^\sharp$ is homotopic to $f$.
\end{lem}

\proof By lemma~\ref{lem-ortho}, the composition of $f$ with the
inclusion of $D$ in the mapping cone of $g$ is nullhomotopic.
Choosing a nullhomotopy and unravelling that gives $f^\sharp\co C\to E$ and
a homotopy from $gf^\sharp$ to $f$. \qed

\begin{defn} We write $\sH\sC'$ for the homotopy category of $\sC'$. By all the above, this is
equivalent to the category $\sC/\sC''$ obtained from $\sC$ by making invertible all morphisms whose mapping
cone belongs to $\sC''$.
\end{defn}

\medskip
In the following lemma, we write $\sC_X$ and $\sC_Y$ etc. instead of $\sC$ to emphasize the dependence of $\sC$ on a space
such as $X$ or $Y$.

\begin{lem} The ``pushforward'' functor $f_*\co \sC_X\to \sC_Y$ determined by $f\co X\to Y$
restricts to a functor $\sC''_X\to \sC''_Y$. If $f$ is an open embedding, it also restricts to
a functor $\sC'_X\to \sC'_Y$. \qed
\end{lem}

\medskip
This completes our discussion of freeness and weak equivalences in $\sC$. We now turn to the
concept of \emph{decomposability}, which is related to excision.

\begin{defn} Let $\sD$ be the smallest full subcategory of $\sC=\sC_X$ with the
following properties.
\begin{itemize}
\item All objects of $\sC$ obtained from
maps $\Delta^k\to X$ (where $k\ge 0$) by the method of
example~\ref{expl-guide} belong to $\sD$.
\item If $C\to D\to E$ is a short exact sequence in $\sC$ and two of the three objects
$C,D,E$ belong to $\sD$, then the third belongs to $\sD$.
\item $\sD\supset \sC''$~, that is, all weakly contractible objects in $\sC$ belong to $0$.
\end{itemize}
When we say that an object of $\sC$ is \emph{decomposable}, we mean that it belongs to $\sD$.
\end{defn}

\begin{defn} $\sD':=\sD\cap \sC'$ and $\sD'':=\sD\cap\sC''$.
\end{defn}

\begin{lem}
\label{lem-retract}
Let $C$ in $\sC$ be obtained from a map $Y\to X$
as in example~\ref{expl-guide}, where $Y$ is a compact ENR. Then $C$ is a retract
in $\sC$ of an object of $\sD$.
\end{lem}

\proof Choose a retraction $r\co Z\to Y$ where $Z$ is a finite
simplicial complex (with right inverse $j\co Y\to Z$, say). Then
$fr\co Z\to X$ determines an object of $\sC$ as in
example~\ref{expl-guide}, and this is clearly in $\sD$. The
object of $\sC$ determined by $f$ is a retract of the object of
$\sD$ determined by $fr$. \qed

\begin{lem} The rule $X\mapsto \sD_X$
is a covariant functor. \qed
\end{lem}

\begin{lem}
\label{lem-approx} Let $V$ and $W$ be open subsets of $X$ such that the
closure of $V$ in $X$ is contained in $W$. Given $D$ in $\sD_X$~,
there are $C$ in $\sD_W$ and a morphism $g\co C\to D$ in $\sD_X$
such that $g_U\co C(U)\to D(U)$ is a homotopy equivalence for any
open $U\subset V$.
\end{lem}

\proof If that claim is true for a particular $D$ and all $V,W$ as
in the statement, then we say that $D$ has property $P$. It is
enough to verify the following.
\begin{itemize}
\item[(i)] Every object obtained from a map $\Delta^k\to X$ by the method of
example~\ref{expl-guide} has property $P$.
\item[(ii)] If $a\co D\to E$ is a weak equivalence in $\sD_X$ and if one of $D,E$ has property $P$,
then the other has property $P$.
\item[(iii)] If $f\co D\to E$ is any morphism in $\sD_X$, and both
$D$ and $E$ have property $P$, then the mapping cone of $f$ has property $P$.
\end{itemize}
The proof of (i) is straightforward using barycentric subdivisions.
Also, one direction of (ii) is trivial: if $D$ has property $P$,
then $E$ has property $P$. For the converse, suppose that $E$ has
property $P$. For fixed $V$ and $W$, choose $g\co F\to E$ with $F$
in $\sD_W$ such that $F(U)\to E(U)$ is a homotopy equivalence
for all $U\subset V$. Without loss of generality, $F$ is free
in every dimension. (Otherwise use lemma~\ref{lem-resolve}.)
Then $F$ belongs to $\sD'_X$ also. By lemma~\ref{lem-cofibrant}, there exists a morphism
$h\co F\to D$ such that the composition $ah\co F\to E$ is homotopic to $g$.
Then $h$ is a morphism which solves our problem. \newline
For the proof of (iii), we fix
$V$ and $W$ and choose an open $W_0$ such that $V\subset W_0\subset
W$, and the closure of $V$ in $X$ is contained in $W_0$ while the
closure of $W_0$ in $X$ is contained in $W$. Then we choose $C$ in
$\sC_W$ and $C\to E$ in $\sC_X$ inducing homotopy equivalences
$C(U)\to E(U)$ for all open $U\subset W_0$. We also choose $B$ in
$\sC_{W_0}$ and $B\to D$ inducing homotopy equivalences $B(U)\to
D(U)$ for all open $U\subset V$. Without loss of generality, $B$
is free in every dimension. Hence there exists $B\to C$ making
the diagram
\[
\CD
B @>>> C \\
@VVV @VVV \\
D @>>> E
\endCD
\]
commutative up to homotopy. Any choice of such a homotopy determines a map
from the mapping cone of $B\to C$ to the mapping cone of $D\to E$
which solves our problem. In particular the mapping cone of
$B\to C$ belongs to $\sC_W$. \qed

\begin{cor}
\label{cor-transverse} Let $X=Y\cup Z$ where $Y,Z$ are open in $X$.
Then for any $D$ in $\sD_X$~, there exists a morphism $C\to D$ in
$\sD_X$ such that $C$ is in $\sD_Y$ and the mapping cone of $C\to D$
is weakly equivalent to an object of $\sD_Z$. \qed
\end{cor}

\proof Choose an open neighborhood $V$ of $X\smin Z$ in $X$ such that
the closure of $V$ (in $X$) is contained in $Y$. Apply lemma~\ref{lem-approx}
with this $V$ and $W=Y$. \qed

\section{Duality in a local setting}

\begin{defn} For an object $C$ in $\sC$ and
open subsets $U,V\subset X$ with $V\subset U$, let $C(U,V)$ be the chain complex
$C(U)/C(V)$  (of free abelian groups).
\end{defn}

\begin{defn}
\label{defn-sheafdual} For objects $C$ and $D$ of $\sC$, let
\[  C\boxtimes D=\holimsub{\twosub{U\subset X~{\rm open},~~K_1,\,K_2\subset X~{\rm closed}}
{K_1 \cap K_2 \subset U}} C(U,U \smin K_1)
\otimes_{\ZZ} D(U,U\smin K_2)~. \]
\end{defn}

The local Poincar\ee duality properties of
topological manifolds do not directly suggest the above definition
of $C\boxtimes D$, but rather an asymmetric definition, as follows.

\begin{defn} For objects $C$ and $D$ of $\sD$, let
\[ C\boxtimes^? D= \holimsub{\twosub{U\subset X~{\rm open}}{K \subset U~{\rm closed}}}
C(U,U\smin K) \otimes_{\ZZ}D(U)~. \]
\end{defn}

\begin{rem} Both $C\boxtimes D$ and $C\boxtimes^?D$ are contractible if either
$C$ or $D$ are in $\sD''$.
\end{rem}

\begin{lem}
\label{lem-forget}
The specialization map $C\boxtimes D\to C\boxtimes^?D$ (obtained by
specializing to $K_2=X$ in the formula for $C\boxtimes D$) induces an isomorphism
in homology.
\end{lem}

\proof We begin with an informal argument. Let $\xi$ be an
$n$-cycle in $C\boxtimes D$. For $K_1, K_2\subset X$ closed,
$U\subset X$ open and $K_1\cap K_2$ contained in $U$,
abbreviate
\[ F(U,K_1,K_2)= C(U,U \smin K_1)
\otimes_{\ZZ} D(U,U\smin K_2). \]
Choose $K^+_1\subset U$
closed, $K^-_1\subset X$ closed so that $K_1=K^+_1\cup K^-_1$ and
$K_2\cap K^-_1=\emptyset$. Then we have a commutative diagram
\[
\xymatrix{ F(U,K^+_1,K_2) \ar[r] & F(U,K^+_1\cap K^-_1,K_2)
& F(U,K^-_1,K_2) \ar[l] \\
F(U,K^+_1,K_2) \ar[r] \ar[u]^= \ar[d]_= & F(U,K^+_1\cap K^-_1,K_2)
\ar[u]^= \ar[d]_=
& F(U\smin K_2,K^-_1,K_2)= 0 \ar[l] \ar[u] \ar[d] \\
F(U,K^+_1,K_2) \ar[r]  & F(U,K^+_1\cap K^-_1,K_2)
& F(U\smin K_2,K^+_1\cap K^-_1,K_2)= 0 \ar[l]  \\
F(U,K^+_1,X) \ar[r] \ar[u] & F(U,K^+_1\cap K^-_1,X) \ar[u] &
F(U\smin K_2,K^+_1\cap K^-_1,X) \ar[l] \ar[u] }
\]
By the sheaf properties, the coordinate of $\xi$ in $F(U,K_1,K_2)$ is
sufficiently determined by the projection of $\xi$ to the homotopy inverse limit of
the top row of the diagram. By diagram chasing, this is sufficiently
determined by the projection of $\xi$ to the homotopy inverse limit of the bottom
row. But that information is stored in the image of $\xi$ in
$C\boxtimes^?D$. \newline
The argument can be formalized as follows. For fixed open $U$, closed $K_1$ and $K_2$
with $K_1\cap K_2$ contained in $U$, we consider the poset $\sP$ of
``decompositions'' $K_1=K^+_1\cup K^-_1$ where $K^+_1\subset U$ and $K^-_1$ are closed,
$K_2\cap K^-_1=\emptyset$. The ordering
is such that $(J^+_1,J^-_1)\le (K^+_1,K^-_1)$ if and only if $J^+_1\subset K^+_1$
and $J^-_1\subset K^-_1$. We need to know that $B\sP$ is contractible. To see this let
$\sQ$ be the poset of all closed neighborhoods of $K_1\cap K_2$ in $K_1\cap U$.
Then $(K^+_1,K^-_1)\mapsto K^+_1$ is a functor $v\co \sP\to \sQ$. Fixing some $J\in \sQ$,
let $\sP_J$ be the poset of all $(K^+_1,K^-_1)\in\sP$ with $K^+_1\subset J$.
This contains as a terminal sub-poset the set of all $(K^+_1,K^-_1)$ with $K^+_1=J$,
and the latter is clearly (anti-)directed. Hence $B\sP_J$ is contractible. This verifies
the hypotheses in Quillen's theorem A for the functor $v$, so that $Bv\co B\sP\to B\sQ$
is a homotopy equivalence. But $\sQ$ is again directed, so $B\sQ$ is contractible.
Therefore $B\sP$ is contractible. Now we can write
\begin{eqnarray*}
C\boxtimes D & = & \holimsub{U,K_1,K_2} F(U,K_1,K_2) \\
&\simeq & \holimsub{U,K_1,K_2} \holimsub{K^+_1,K^-_1} F(U,K_1,K_2) \\
&\cong& \holimsub{U,K^+_1,K^-_1,K_2} F(U,K^+_1\cup K^-_1,K_2) \\
&\simeq& \holimsub{U,K^+_1,K^-_1,K_2} \holim
\left[F(U,K^+_1,K_2) \to F(U,K^+_1\cap K^-_1,K_2)
\leftarrow F(U,K^-_1,K_2)\right].
\end{eqnarray*}
(The usual conventions apply: $K_1,K_2$ closed in $X$, with
$K_1\cap K_2$ contained in the open set $U$, and $K_1=K^-_1\cup K^+_1$ is a
decomposition of the type we have just discussed.)
Using that, we obtain from the rectangular
twelve term diagram above a map $g\co H_*(C\boxtimes^?D)\to H_*(C\boxtimes D)$.
Indeed an $n$-cycle in $C\boxtimes^?D$ determines an $n$-cycle in
\[
\holimsub{U,K^+_1,K^-_1,K_2} \holim\left[
F(U,K^+_1,X) \to F(U,K^+_1\cap K^-_1,X) \leftarrow F(U\smin K_2,K^+_1\cap K^-_1,X)\right]
\]
and we use the vertical arrows in the twelve-term diagram to obtain
an $n$-cycle in
\[
\holimsub{U,K^+_1,K^-_1,K_2} \holim
\left[F(U,K^+_1,K_2) \to F(U,K^+_1\cap K^-_1,K_1)
\leftarrow F(U,K^-_1,K_2)\right]
\]
whose homology class is well defined. By construction,
$g\co H_*(C\boxtimes^?D)\to H_*(C\boxtimes D)$ is left inverse to the projection-induced
map $H_*(C\boxtimes D)\to H_*(C\boxtimes^?D)$. But it is clearly also right inverse
(specialize to $K_2=X$ and then $K^-_1=\emptyset$ in the diagrams above). \qed

\begin{expl}
\label{expl-can1} Let $X$ be a compact ENR. For open $U$ in $X$,
let $C(U)$ be the singular chain complex of $U$, regarded as a subcomplex of $C$. By
lemma~\ref{lem-guide}, applied to the identity $X\to X$, this
functor $U\mapsto C(U)$ satisfies the sheaf type condition and the two
finiteness conditions. We construct a ``canonical'' map
\[  \nabla\co C(X) \lra C\boxtimes C~.  \]
To start with we have the following diagram:
\[
\CD
C(X) @>>>  {\holimsub{U,K} C(X,X\smin K)} @<\simeq<<  {\holimsub{U,K} C(U,U\smin K)}
\endCD
\]
where $U$ and $K$ are open in $X$ and closed in $X$, respectively, with $K\subset U$. The first
arrow is induced by the quotient maps $C(X)\to C(X,X\smin K)$ and the second map is induced
by the inclusions $C(U,U\smin K)\to C(X,X\smin K)$ which are chain homotopy equivalences by the
sheaf property. Inversion of the second arrow in the diagram gives us a map
\[ \gamma\co C(X) \lra \holimsub{U,K} C(U,U\smin K)\]
well defined up to contractible choice (in particular, well defined up to chain homotopy). \newline
Next we make use of a certain chain map
\[ \zeta\co \holimsub{U,K} C(U,U\smin K)\la C\boxtimes C~. \]
This is determined by the compositions
\[
\CD
C(U,U\smin(K_1\cap K_2)) @>>>
\textup{sing. chain cx. of } (U\times U,U\times U\smin(K_1\times U\,\cup\,U\times K_2)) \\
@.     @VVV \\
 @. C(U,U\smin K_1)\otimes C(U,U\smin K_2)
\endCD
\]
(for closed $K_1,K_2\subset X$ with intersection $K_1\cap K_2\subset U$), where
the first arrow is induced by the diagonal map
and the second arrow is an Eilenberg-Zilber map. Now we define
\[  \nabla=\zeta\gamma\co C(X) \lra C\boxtimes C~. \]
This is a refinement of the standard Eilenberg-Zilber-Alexander-Whitney diagonal chain map
$C(X)\to C(X)\otimes C(X)$, which we can recover by composing with the projection (alias specialization)
from $C\boxtimes C$ to $C(X)\otimes C(X)$. \newline
If $X$ is an oriented closed topological $n$-manifold and $\omega\in C(X)$ is an $n$-cycle representing
a fundamental class for $X$, then $\nabla(\omega)$ is an $n$-cycle in $C\boxtimes C$ which, as we shall
see in proposition~\ref{prop-nondeg} below, is ``nondegenerate''.
This reflects the fact that not only $X$, but also each open subset of $X$ satisfies
a form of Poincar\ee duality. \end{expl}

\begin{expl}
\label{expl-can2}
Let $(X,Y)$ be a pair of compact ENRs. For open $U$ in $X$, let $C(U)$ be the singular chain complex of $U$
and let $D(U)$ be the singular chain complex of $U\cap Y$. By lemma~\ref{lem-guide}, both $C$ and $D$ satisfy the
sheaf type condition and the two finiteness conditions. A straightforward generalization of the
previous example gives
\[ \nabla\co C(X)/D(X) \lra (C\boxtimes C)\big/(D\boxtimes D)~. \]
If $X$ is an oriented compact topological $n$-manifold with boundary $Y$,
and $\omega\in C(X)/D(X)$ is an $n$-cycle representing
a fundamental class for the pair $(X,Y)$, then $\nabla(\omega)$ is an $n$-cycle in
$(C\boxtimes C)\big/(D\boxtimes D)$ whose image in $C\boxtimes(C/D)$
is nondegenerate (proposition~\ref{prop-nondeg} below).
\end{expl}

\begin{prop}
\label{prop-nondeg}
Let $C$ and $D$ be objects of $\sC$. Let $[\varphi]\in H_n(C\boxtimes D)$
and suppose that, for every open $U\subset X$ and every $j\in \ZZ$, the map
\[  \colimsub{\textup{cpct }K\subset U} H^{n-j}C(U,U\smin K)\la H_jD(U)~,   \]
slant product with the coordinate of $\varphi$ in $C(U,U\smin K)\otimes D(U)$,
is an isomorphism. Then $[\varphi]$ is nondegenerate in the
following sense: for any $E$ in $\sD$, the map $f\mapsto f_*[\varphi]$ is an isomorphism
from $[D,E]$, the morphism set in $\sC/\sC''\cong \sH\sC'$, to $H_n(C\boxtimes E)$.
\end{prop}

\proof The idea is very simple. Given $[\psi]\in H_n(C\boxtimes E)$, the slant product
with $[\psi]$ gives us homomorphisms
\[ \colimsub{\textup{closed }K\subset U} H^{n-j}C(U,U\smin K)\la H_jE(U) \]
for $j\in \ZZ$. By our assumption on $[\varphi]$, we have
$\colim_K H^{n-j}C(U,U\smin K)\cong H_jD(U)$
and so we get homomorphisms
\[  H_jD(U) \la H_jE(U) \]
for all $j\in \ZZ$, naturally in $U$. It remains ``only'' to construct a morphism
$f_{\psi}\co D\to E$ in $\sD$ or in $\sH\sD/\sH\sD''\cong \sH\sD'$ inducing these homomorphisms
of homology groups. \newline
By lemma~\ref{lem-forget} or otherwise, we may represent the class $[\psi]$ by an $n$-cycle
\[ \psi \in C\boxtimes^?E=\holimsub{\twosub{U\subset X~{\rm open}}{K \subset U~{\rm closed}}}
C(U,U\smin K) \otimes_{\ZZ}E(U). \]
We shall show first of all that this $n$-cycle determines a chain map
\[
\psi^{\textup{ad}}\co
\hocolimsub{K\subset U} C(X,X\smin K)^{n-*} \lra \holimsub{V\supset U} E(V) \]
natural in the variable $U$. In fact the source of this map is clearly a covariant
functor of the variable $U$, by \emph{extension of the indexing poset}.
The target is homotopy equivalent to $E(U)$, and we regard it as a covariant functor of $U$
by \emph{restriction of the indexing poset}. More precisely, if $U_1\subset U_2$,
then the poset of open subsets of $X$ containing $U_2$ is contained in the subposet
of open subsets of $X$ containing $U_1$. Corresponding to that inclusion of posets we have
a projection map
\[  \holimsub{V\supset U_1} E(V)  \la \holimsub{V\supset U_2} E(V) \]
and that is what we use. \newline
To give a precise description of $\psi^{\textup{ad}}$ now, we fix a string
$K_0\subset K_1\subset \dots\subset K_r$
of compact subsets of $X$, and a string of open subsets
$U_s\subset U_{s-1}\subset \cdots\subset U_0$
with $K_r\subset U_s$. Let $\ell\Delta^r$ and $\ell\Delta^s$ be the
cellular chain complexes of $\Delta^r$ and $\Delta^s$ (the $\ell$ is for
\emph{linearization}). We can define $\psi^{\textup{ad}}$ by associating to every choice
of two such strings a chain map
\[ \ell\Delta^r\otimes_\ZZ C(X,X\smin K_0)^{n-*} \la \hom_\ZZ(\ell\Delta^s,E(U_0)), \]
or equivalently, a chain map
\[ \ell\Delta^s\otimes_\ZZ\ell\Delta^r \la C(X,X\smin K_0)\otimes_\ZZ E(U_0) \]
of degree $n$.
(As the input strings vary, these chain maps are subject to obvious compatibility conditions.)
But in fact the coordinate of $\psi\in C\boxtimes^? E$ corresponding to the strings
$(K_0,\dots,K_r)$ and $(U_s,\dots,U_0)$ is a chain map
\[ \ell\Delta^s\otimes\ell\Delta^r \la C(U_0,U_0\smin K_0)\otimes E(U_0)  \]
of degree $n$. We need only compose with the inclusion-induced maps
\[ E(U_0)\otimes D(U,U\smin K_r) \la E(U_0)\otimes D(X,X\smin K_0) \]
to get the data we need. \newline
Now abbreviate
\begin{eqnarray*}
PC^{(n-*)}(U) & = & \hocolimsub{K\subset U} C(X,X\smin K)^{n-*}, \\
PE(U) & = & \holimsub{V\supset U} E(V), \\
\quad PD(U) & = & \holimsub{V\supset U} D(V)
\end{eqnarray*}
(writing $P$ for \emph{provisional}).
Then we have a homotopy commutative diagram of natural transformations of functors
on $\sO(X)$,
\[
\CD
D @<<< C^{(n-*)} @>>>  E \\
@VVV  @VVV     @VVV  \\
PD @<\varphi^{\textup{ad}}<< PC^{(n-*)} @>\psi^{\textup{ad}}>> PE
\endCD
\]
as follows. The maps $E\to PE$ and $D\to PD$ are obvious (diagonal) constructions.
We note that $E(U)\to PE(U)$ and $D(U)\to PD(U)$ are homology equivalences
for all $U$. The object $C^{(n-*)}$ is chosen in $\sC'$ and the
map from it to $PC^{(n-*)}$ induces homology equivalences
$C^{(n-*)}(U)\to PC^{(n-*)}(U)$ for all $U$ by construction. (This uses
lemma~\ref{lem-resolve}.) The arrows
in the top row are then constructed to make the diagram homotopy
commutative. (This uses lemma~\ref{lem-cofibrant}.)
The horizontal arrows in the left-hand square of the
diagram are homology equivalences (for every choice of input $U$) and
consequently the arrow from $C^{(n-*)}$ to $D$ is an isomorphism
in $\sC/\sC''$. We choose an inverse for it and compose with
$C^{(n-*)}\to E$ to obtain the desired element $[f_\psi]\in [D,E]$.
Finally, it is just a matter of inspection to see that the rule
$[\psi]\mapsto [f_\psi]$ is inverse to the map from $[D,E]$ to
$H_n(C\boxtimes E)$ given by slant product with $[\varphi]$.
\qed

\begin{expl}
\label{expl-can3} Recall that the natural Eilenberg-Zilber map
\[ \textup{sg ch cx of }(Y\times Z) \la
(\textup{sg ch cx of }Y)\otimes (\textup{sg ch cx of }Z) \]
(for arbitrary spaces $Y$ and $Z$) admits a refinement to a natural
equivariant chain map
\[ W\otimes(\textup{sg ch cx of }(Y\times Z))
\la (\textup{sg ch cx of }Y)\otimes (\textup{sg ch cx of }Z) \]
where $W$ is a free resolution of $\ZZ$ as a trivial module over
the group ring $\ZZ[\ZZ/2]$. Here \emph{equivariance} means that
\[
\CD
W\otimes(\textup{sg ch cx of }Y\times Z)
@>>> (\textup{sg ch cx of }Y)\otimes (\textup{sg ch cx of }Z) \\
@V T\otimes \textup{perm.} V\cong V @V \textup{perm.} V\cong V \\
W\otimes(\textup{sg ch cx of }Z\times Y)
@>>> (\textup{sg ch cx of }Z)\otimes (\textup{sg ch cx of }Y)
\endCD
\]
commutes for all $Y$ and $Z$, where $T$ is the generator of $\ZZ/2$ acting
on $W$ and ``perm.'' stands for a permutation of the factors $Y$ and $Z$. \newline
Applying this in the situation of example~\ref{expl-can1}, we deduce immediately
that
\[ \nabla\co C(X) \la C\boxtimes C \]
has a refinement to
\[ \nabla^{h\ZZ/2}\co C(X) \la (C\boxtimes C)^{h\ZZ/2}. \]
Suppose now that $X$ is a closed topological manifold, and $\omega\in C(X)$
is a fundamental cycle. Then with the nondegeneracy property
which we have already established, we can say informally that
$(C,\nabla^{h\ZZ/2}(\omega))$ is an $n$-dimensional SAPC in $\sC$ (slightly against
the rules, since we have not established that every object in $\sC$ has a dual). \newline
Similarly, in the situation and notation of example~\ref{expl-can2},
we obtain $\nabla^{h\ZZ/2}(\omega)$, an $n$-cycle in
$(C\boxtimes C)^{h\ZZ/2}/(D\boxtimes D)^{h\ZZ/2}$. With the nondegeneracy property
which we have already established, this allows us to say informally that
$((C,D),\nabla^{h\ZZ/2}(\omega))$ is an SAP pair.
\end{expl}

\medskip
There are ``economy'' versions of
$C\boxtimes D$ other than $C\boxtimes^?D$. Suppose that $\sQ$ is a collection of closed subsets
of $X$ with the following properties.
\begin{itemize}
\item $\sQ$ is closed under finite intersections (i.e., for $Q_1$ and $Q_2$ in $\sQ$,
we have $Q_1\cap Q_2\in \sQ$)~;
\item for every compact subset $K$ of $X$ and
open $U\subset X$ containing $K$, there exist $r\ge 0$ and $Q_1,\dots,Q_r\in \sQ$ such
that
\[   K\subset \bigcup_{i=1}^r Q_r \subset U~. \]
\end{itemize}
For $C$ and $D$ in $\sD$ we define
\[  C\boxtimes^{??}D=\holimsub{\twosub{\twosub{U\in\sO(X)}{K_1,K_2\in \sQ}}{K_1\cap K_2\subset U}}
C(U,U\smin K_1)\otimes_\ZZ D(U,U\smin K_2)~. \]
This depends on $\sQ$, not just on $C$ and $D$, but in practice it will be
clear what $\sQ$ is.

\begin{lem}
\label{lem-prune}
The specialization map $C\boxtimes D\la C\boxtimes^{??}D$ is a chain homotopy equivalence.
\end{lem}

\proof Step 1: We assume that $X$ is compact.
Let $\sQ'$ be the collection of all subsets of $X$ which are finite unions of subsets
$K\in \sQ$. The specialization map $C\boxtimes D\la C\boxtimes^{??}D$ is a composition of two
specialization maps
\[
\xymatrix{
{\holimsub{\twosub{\twosub{U\in\sO(X)}{K_1,K_2\textup{ closed}}}{K_1\cap K_2\subset U}}
C(U,U\smin K_1)\otimes_\ZZ D(U,U\smin K_2)} \ar[d]_f \\
{\rule{0mm}{6mm}\holimsub{\twosub{\twosub{U\in\sO(X)}{K_1,K_2\in \sQ'}}{K_1\cap K_2\subset U}}
C(U,U\smin K_1)\otimes_\ZZ D(U,U\smin K_2)} \ar[d]_g \\
{\rule{0mm}{6mm}\holimsub{\twosub{\twosub{U\in\sO(X)}{K_1,K_2\in \sQ}}{K_1\cap K_2\subset U}}
C(U,U\smin K_1)\otimes_\ZZ D(U,U\smin K_2).}
}
\]
We are going to show that both of these are homotopy equivalences. For the specialization map $f$,
it suffices to observe that, by our assumptions on $\sQ$, the triples
$(U,K_1,K_2)$ with $K_1,K_2\in \sQ'$ and
$K_1\cap K_2\subset U$ form an \emph{initial} sub-poset in the poset of all triples
$(U,K_1,K_2)$ with $K_1,K_2$ closed and $K_1\cap K_2\subset U$. This refers to the usual ordering,
\[  (U,K_1,K_2) \le (V,L_1,L_2)\qquad\Leftrightarrow\qquad U\subset V \textup{ and }
K_1\supset L_1~, K_2\supset L_2~. \]
For the specialization map $g$, it suffices by \cite[9.7]{DwyerKan} to show that
for open $U\subset X$ and $L_1,L_2$ in $\sQ'$ with $L_1\cap L_2\subset U$, the canonical map
\[ C(U,U\smin L_1)\otimes_\ZZ D(U,U\smin L_2) \lra
\holimsub{\twosub{K_1,K_2\in \sQ}{K_1\subset L_1\,,\,K_2\subset L_2}}
C(U,U\smin K_1)\otimes_\ZZ D(U,U\smin K_2)
\]
is a chain homotopy equivalence. But this is true by
lemma~\ref{lem-descent}. (\emph{Some details}: The target of this map can also be described,
up to homotopy equivalence, as a double homotopy limit
\[  \holimsub{\sR}\holimsub{\twosub{K_1,K_2\in \sR}{K_1\subset L_1\,,\,K_2\subset L_2}}
C(U,U\smin K_1)\otimes D(U,U\smin K_2) \]
where $\sR$ runs through the \emph{finite} subsets of $\sQ$ which are ``large enough''.
By \emph{large enough} we mean that there exist $K_{11}, K_{12}, \dots, K_{1r}$ and
$K_{21}, K_{22}, \dots, K_{2s}$ in $\sR$ such that
\[  \bigcup_{i=1}^r K_{1i}=L_1~,~~\bigcup_{j=1}^r K_{2j}=L_2~. \]
For fixed $\sR$, we have an
Alexander-Whitney type homotopy equivalence
\[
\begin{array}{cl}
& \rule{0mm}{5mm}\holimsub{\twosub{K_1,K_2\in \sR}{K_1\subset L_1\,,\,K_2\subset L_2}}
C(U,U\smin K_1)\otimes D(U,U\smin K_2) \\
\simeq & \rule{0mm}{8mm}\bigg(\holimsub{\twosub{K_1\in \sR}{K_1\subset L_1}}
C(U,U\smin K_1)\bigg)   \otimes \bigg(\holimsub{\twosub{K_2\in \sR}{K_2\subset L_2}}
D(U,U\smin K_2)\bigg)
\end{array}
\]
natural in $\sR$. By lemma~\ref{lem-descent} and because $\sR$ is large enough, the projections
\[
\begin{array}{ccc}
C(U,U\smin L_1) & \lra & \holimsub{\twosub{K_1\in \sR}{K_1\subset L_1}} C(U,U\smin K_1)~, \\
D(U,U\smin L_2) & \lra & \holimsub{\twosub{K_2\in \sR}{K_2\subset L_2}} D(U,U\smin K_2)
\end{array}
\]
are homotopy equivalences. Putting these facts together, we see that
\[ \holimsub{\twosub{K_1,K_2\in \sQ}{K_1\subset L_1\,,\,K_2\subset L_2}}
C(U,U\smin K_1)\otimes_\ZZ D(U,U\smin K_2) \]
is homotopy equivalent to the homotopy inverse limit of a constant functor,
\[ \sR~~\mapsto~~C(U,U\smin K_1)\otimes_\ZZ D(U,U\smin K_2). \]
Since the poset of all
$\sR$ is directed, the homotopy inverse limit is homotopy equivalent to the
unique value of that functor.) \newline
Step 2: $X$ is arbitrary (but still locally compact Hausdorff and separable). Choose a compact $Y\subset X$
which belongs to $\sQ$ and is a neighborhood for the support of $C$ and for the support of $D$. Let
\[ \sQ^Y=\{ Q\in \sQ~|~ Q\subset Y \}~, \]
a collection of compact subsets of $Y$ which is closed under finite intersections.
For open $U\subset Y$ put
\[ C^Y(U)=C(U\cup(X\smin Y))~,\quad D^Y(U)=D(U\cup(X\smin Y))~. \]
Now we have a commutative diagram of specialization maps
\[
\xymatrix{
C\boxtimes D \ar[r] \ar[d] & C\boxtimes^{??}D  \ar[d] \\
C^Y\boxtimes D^Y \ar[r] & C^Y\boxtimes^{??}D^Y
}
\]
using $\sQ^Y$ to define the lower row. By step 1, the lower horizontal arrow is a homotopy
equivalence. It is therefore enough to show that the two vertical arrows are homotopy
equivalences. This follows easily from the fact that the inclusion
of posets $\iota:\sW\to \sO(X)$ has a left adjoint, where $\sO(X)$ consists of all open
subsets of $X$ and $\sW$ consists of all open subsets of $X$ containing $X\smin Y$.
The left adjoint is given by $U\mapsto \lambda(U)=U\cup(X\smin Y)$.
Note also that the inclusion-induced maps $C(U)\cong C(\lambda(U))$ and
$D(U)\cong D(\lambda(U))$ are isomorphisms. Thus, $\lambda$ induces maps
$C^Y\boxtimes D^Y\to C\boxtimes D$ and $C^Y\boxtimes^{??}D^Y\to C\boxtimes^{??}D$ which
are homotopy inverses for the vertical arrows in our square. The homotopies are induced
by natural transformations, the ``unit'' and the ``counit'' of the adjunction
of $\iota$ and $\lambda$.
\qed

\section{Products}
Let $X$ and $Y$ be locally compact Hausdorff and separable spaces.

\begin{defn}
\label{defn-product} For $C$ in $\sC_X$ and $D$ in $\sC_Y$,
the \emph{tensor product} $C\otimes D$ of $C$ and $D$ is the ordinary tensor product
of chain complexes $C\otimes_\ZZ D$~, with the system of subcomplexes defined by
\[ (C\otimes D)(W)= \sum_{\twosub{U,V}{U\times V\subset W}} C(U)\otimes_\ZZ D(V) \]
for $W\in \sO(X\times Y)$.
\end{defn}

We are aiming to show that $C\otimes D$ is in $\sC_{X\times Y}$. This is
surprisingly hard. We begin with two lemmas.

\begin{lem}
\label{lem-hocoandsum}
For $i\in\{1,2,\dots,k\}$ let $U_i$ be open in $X$ and let $V_i$ be
open in $Y$. For nonempty $S\subset\{1,2,\dots,k\}$ put
$U_S=\bigcap_{\lambda\in S}U_\lambda$ and
$V_S=\bigcap_{\lambda\in S}V_{\lambda}$. The following map (induced by
obvious inclusions) is a homotopy equivalence:
\[ \hocolimsub{\emptyset\ne S\subset\{1,2,\dots,k\}}
C(U_S)\otimes_\ZZ D(V_S) \la \sum_{i=1}^k C(U_i)\times D(V_i)~, \]
where the sum $\sum_{i=1}^k$ is taken inside $C(X)\otimes_\ZZ D(Y)$.
\end{lem}

\proof We proceed by induction on $k$.
The square of inclusion maps
\[
\CD
\sum_{i=1}^{k-1} C(U_i\cap U_k)\otimes D(V_i\cap V_k)  @>>> C(U_k)\otimes D(V_k) \\
@VVV   @VVV  \\
\sum_{i=1}^{k-1} C(U_i)\otimes D(V_i)  @>>> \sum_{i=1}^k C(U_i)\otimes D(V_i)
\endCD\quad (*)
\]
is a homotopy pushout square. Indeed, it is a pushout square in which
the horizontal arrows (in fact, also the vertical arrows) are cofibrations, i.e.,
split injective as maps of graded abelian groups. The pushout property can be
verified in terms of bases: each of the four terms in the square is the graded free
abelian group generated by a certain graded set. \newline
Next, for nonempty $S\in \{1,\dots,k-1\}$, write $E(S)=
C(U_S)\times D(V_S)$. Then it is clear that
\[
\CD
\hocolimsub{S} E(S\cup k)
@>>> \hocolimsub{S} E(k) \\
@VVV @VVV \\
\hocolimsub{S} E(S)
@>>> \hocolimsub{S} \hocolim(E(S)\leftarrow E(S\cup k)\to E(k))~,
\endCD (**)
\]
where $S$ refers to nonempty subsets of $\{1,\dots,k-1\}$~,
commutes up to a preferred homotopy $h$ and, as such, is a homotopy pushout square.
The square $(**)$ maps to $(*)$ by a forgetful map which also takes the
homotopy $h$ to zero. By inductive hypothesis, three of the four
arrows which constitute this map $(**)\to (*)$ between squares are homotopy equivalences, and therefore
all are homotopy equivalences. Now it only remains to show that the canonical map
\[
\CD
\hocolimsub{\emptyset\ne S\subset\{1,2,\dots,k-1\}} \hocolim(E(S)\leftarrow E(S\cup k)\to E(k)) \\
@VVV \\
\hocolimsub{\emptyset\ne S\subset\{1,2,\dots,k\}} E(S)
\endCD
\]
is a homotopy equivalence. We have dealt with this kind of task before, in the proof of
lemma~\ref{lem-descent} and lemma~\ref{lem-codescent}, and it can be dealt with in the same way here.
\qed

\begin{lem}
\label{lem-sumapprox}
In the situation of lemma~\ref{lem-hocoandsum}, let $W$ be the union of
the sets $U_i\times V_i$ for $i=1,\dots,k$. Then the inclusion
\[   \sum_{i=1}^k C(U_i)\otimes_\ZZ D(V_i) \la (C\otimes D)(W) \]
is a homotopy equivalence.
\end{lem}

\proof Step 1: We assume to begin with that $W$ itself has the form $U\times V$ for some
$U$ open in $X$ and some $V$ open in $Y$. It is easy to construct finite open coverings
$\{U'_\lambda~|~\lambda\in \Lambda\}$ of $U$, and $\{V'_\gamma~|~\gamma\in \Gamma\}$ of $V$,
such that every open set $U_i\times V_i$ is a union of (some) of the sets $U'\lambda\times V'_\mu$.
Now the composition of inclusions
\[
\CD
{\ds\sum_{(\lambda,\mu)} C(U'_\lambda)\otimes D(V'_\mu)} \\
@VVV \\
{\ds \sum_{i} C(U_i)\otimes D(V_i)} \\
@VVV \\
C(U\times V)
\endCD
\]
is a homotopy equivalence, because the source can be written as
\[  \left(\sum_\lambda C(U'_\lambda)\right)\otimes_\ZZ \left(\sum_\mu C(V'_\mu)\right) \]
and we are assuming the sheaf type condition for $C$ and $D$. Therefore it remains
only to show that the first of these inclusions admits a homotopy left inverse.
Using lemma~\ref{lem-hocoandsum}, we may replace $\sum_i C(U_i)\otimes D(V_i)$
by
\[
\hocolimsub{\emptyset\ne T \subset\{1,\dots,k\}}
C(U_T)\otimes D(V_T).
\]
In that expression, $C(U_T)$ can be replaced by the sum of the $C(U'_R)$
for $U'_R\subset U_T$, where $R$ is a nonempty subset of $\Lambda$. Similarly
$D(V_T)$ can be replaced by the sum of the $D(V'_S)$
for $V'_S\subset V_T$, where $S$ is a nonempty subset of $\Gamma$. (Here we use the
sheaf type conditions for $C$ and $D$ again.) After these modifications, we have an
obvious projection map from that homotopy colimit to
\[ \sum_{(\lambda,\mu)} C(U'_\lambda)\otimes D(V'_\mu)~. \]
This provides the required left homotopy inverse. \newline
Now we look at the general case. Let $W=\bigcup_{i=1}^k U_i\times V_i$ as given. Let $U_{k+1}\subset X$
and $V_{k+1}\subset Y$ be open and suppose $U_{k+1}\subset V_{k+1}\subset W$. It is enough to show that
the inclusion
\[  \sum_{i=1}^k C(U_i)\otimes D(V_i) \la \sum_{i=1}^{k+1} C(U_i)\otimes D(V_i) \]
is a homotopy equivalence. (For then we can repeat the process by adding on
as many terms $C(U')\otimes D(V')$ with $U'\times V'\subset W$ as we like, and thereby
approximate $(C\otimes D)(W)$.) There is a pushout square
\[
\CD
\sum_{i=1}^k C(U_i\cap U_{k+1})\otimes D(V_i\cap V_{k+1})  @>>> C(U_{k+1})\otimes D(V_{k+1}) \\
@VVV  @VVV  \\
\sum_{i=1}^k C(U_i)\otimes D(V_i) @>>> \sum_{i=1}^{k+1} C(U_i)\otimes D(V_i)
\endCD
\]
in which the horizontal arrows are cofibrations. It follows that the square
is a homotopy pushout square. By step 1, the upper horizontal arrow is a homotopy equivalence.
Therefore the lower horizontal arrow is a homotopy equivalence.
\qed

\begin{lem}
\label{lem-productsheaf}
If $C$ belongs to $\sC_X$ and $D$ belongs to $\sC_Y$~,
then $C\otimes D$ belongs to $\sC_{X\times Y}$.
\end{lem}

\proof For the sheaf condition, suppose that $W\subset X\times Y$ is open and $W=\bigcup W_\alpha$. Then
$W$ is the union of all open sets $U_i\times V_i$ which are contained in some $W_\alpha$,
and therefore the inclusion
\[  \sum_\alpha (C\otimes D)(W_\alpha) = \sum_i C(U_i)\otimes D(U_i) \la (C\otimes D)(W) \]
is a homotopy equivalence by lemma~\ref{lem-sumapprox} and passage to the direct limit. \newline
For finiteness condition (ii), suppose that $C$ has support in a compact
subset $K$ of $X$ and $D$ has support in a compact subset $L$ of $Y$. Then it is clear that
$C\otimes D$ has
support in $K\times L\subset X\times Y$. \newline
This leaves finiteness condition (i) to be established. We recall what it requires. We have to
find an integer $c\ge 0$ such that, for open $W\subset W'$ in $X\times Y$, where $W$ has compact
closure in $W'$, the inclusion $(C\otimes D)(W)\to (C\otimes D)(W')$ factors through a chain
complex of f.g. free abelian groups whose $i$-th chain group is zero whenever $|i|>c$. Let
$a,b\in \ZZ$ be the corresponding integers for $C$ and $D$. Let us provisionally say that an
open set $W$ in $X\times Y$ is \emph{good} if
\begin{itemize}
\item[] $W$ has compact closure in $X$ and, for every open $W''$ in $X$ containing
the closure $\bar W$ of $W$, there exists another open $W'$ with compact
closure in $X$ such that
\[  W\subset \bar W\subset W'\subset \bar W'\subset W \]
and the inclusion map $(C\otimes D)(W')\to (C\otimes D)(W'')$ factors
through a bounded chain complex of f.g. free abelian groups.
\end{itemize}
Then it is easy to verify:
\begin{itemize}
\item any $W$ of the form $W=U\times V$, with compact closure in $X\times Y$, is good;
\item if $W_1$, $W_2$ and $W_1\cap W_2$ are good open subsets of $X\times Y$, then
$W_1\cup W_2$ is good.
\end{itemize}
It follows that if $W$ is any finite union of subsets of the form
$U\times V$, where $U$ and $V$ are open in $X$ and $Y$,
respectively, with compact closures, then $W$ is good. While that does not prove
all we need, it will now be sufficient to show that, for any $W$ open in
$X\times Y$, the chain complex $(C\otimes D)(W)$ is homotopy equivalent to
a chain complex of free abelian groups concentrated in degrees between $-ab$ and
$(a+2)(b+2)+1$. To that end,
we note that for any open $U$ in $X$, the chain complex $C(U)$ is homotopy
equivalent to a chain complex of free abelian groups concentrated in degrees
between $-a$ and $a+1$, as it is homotopy equivalent to
a sequential homotopy colimit of chain complexes concentrated in degrees between $-a$ and $a$.
Consequently, for open $U',U''$ in $X$ with $U'\subset U''$, the chain complex
$C(U',U'')=C(U')/C(U'')$ is homotopy equivalent to a chain complex of free abelian
groups concentrated in degrees between $-a$ and $a+2$. Similar remarks apply to $D$
in place if $C$. If $W$ is open in $X\times Y$ and is a finite union
of subsets of the form $U_\alpha\times V_\alpha$, then we may replace $(C\otimes D)(W)$
by the homotopy equivalent subcomplex
\[  \sum_\alpha C(U_\alpha)\otimes D(V_\alpha). \]
This admits a finite filtration
by subcomplexes such that the subquotients of the filtration have the form
\[ C(U',U'')\otimes_\ZZ D(V',V'') \]
for some open $U',U''$ in $X$ and $V',V''$ in $Y$, with $U'\subset U''$ and $V'\subset V''$.
It is therefore homotopy equivalent to a chain complex of free abelian groups
concentrated in degrees between $-ab$ and $(a+2)(b+2)$. Finally, an arbitrary
open $W$ in $X\times Y$ can be written as a monotone union of subsets $W_i$
where each $W_i$ is a finite union of open subsets of the form $U_\alpha\times V_\alpha$.
Since each $(C\otimes D)(W_i)$ is homotopy equivalent to a chain complex of free
abelian groups concentrated in degrees between $-ab$ and $(a+2)(b+2)$, it follows that
$(C\otimes D)(W)$ is homotopy equivalent to a chain complex of free
abelian groups concentrated in degrees between $-ab$ and $(a+2)(b+2)+1$.
\qed

\begin{cor}
\label{cor-productsheaf} If $C$ belongs to $\sD_X$ and $D$ belongs
to $\sD_Y$~, then $C\otimes D$ belongs to $\sD_{X\times Y}$.
\end{cor}

\proof It is clear from the definition that, for fixed $C$, the functor
$D\mapsto C\otimes D$ from $\sC_X$ to $\sC_{X\times Y}$ respects weak equivalences
and short exact sequences. The same is true if we fix $D$ and allow $C$ to vary.
Therefore it is enough to prove the claim in the special case where $C$ and $D$ are
obtained from maps $f\co \Delta^k\to X$ and $g\co \Delta^\ell\to Y$, so that $C(U)$
and $D(V)$ are the singular chain complexes of $f^{-1}(U)$ and $g^{-1}(V)$, respectively,
for $U$ open in $X$ and $V$ open in $Y$. Let $E$ in $\sC_{X\times Y}$ be the object
obtained from $f\times g\co \Delta^k\times\Delta^\ell\to X\times Y$, so that $E(W)$
is the singular chain complex of $(f\times g)^{-1}(W)$, for $W$ open in $X\times W$.
It is easy to show that $E$ belongs to $\sD_{X\times Y}$. There is an easy Eilenberg-Zilber
type map
\[  C\otimes D \lra E \]
in $\sC_{X\times Y}$. For open $W$ in $X\times Y$ of the form $W=U\times V$ with
$U$ open in $X$ and $V$ open in $Y$, this specializes to a map of chain complexes
\[  (C\otimes D)(W) \lra E(W) \]
which is a chain homotopy equivalence by the very Eilenberg-Zilber theorem. It follows
then from the sheaf properties that $(C\otimes D)\to E(W)$ is always a chain homotopy
equivalence, for arbitrary open $W$ in $X\times Y$.
Consequently $C\otimes D$ belongs to $\sD_{X\times Y}$. \qed

\begin{prop}
\label{prop-boxes}
Let $C$ and $E$ be objects of $\sC_X$. Let $D$ and $F$ be objects of $\sC_Y$. Then the
following specialization map is a homotopy equivalence (and a ``fibration'', i.e.,
split surjective as a map of graded abelian groups):
\[
\xymatrix@C=5pt@R=35pt{
(C\otimes D)\boxtimes(E\otimes F) &  = & {\rule{0mm}{5mm}\holimsub{W,P_1,P_2}
(C\otimes D)(W,W\smin P_1)\otimes (E\otimes F)(W,W\smin P_2)}  \ar[d] \\
 &  &
{\rule{0mm}{6mm}\holimsub{\twosub{\twosub{W=U\times V}{P_1=J_1\times K_1}}{P_2=J_2\times K_2}}\!\!\!\!
(C\otimes D)(W,W\smin P_1) \otimes (E\otimes F)(W,W\smin P_2).}
}
\]
\end{prop}

\proof
As in the proof of lemma~\ref{lem-prune}, there is no loss of generality
in assuming that $X$ and $Y$ are both compact. In that case we may replace the target
of our specialization map by
\[ \holimsub{\twosub{\twosub{W}{P_1=J_1\times K_1}}{P_2=J_2\times K_2}}
(C\otimes D)(W,W\smin P_1) \otimes (E\otimes F)(W,W\smin P_2)~, \]
dropping the condition $W=U\times V$. (This makes no difference to the homotopy
type because, in the poset of triples $(W,P_1,P_2)$ where $P_1=J_1\times K_1$
and $P_2=J_2\times K_2$, those triples which have $W=U\times V$ for some $U$ and
$V$ form an ``initial'' sub-poset.) Now our map has the form
\[ (C\otimes D)\boxtimes(E\otimes F) \la (C\otimes D)\boxtimes^{??}(E\otimes F) \]
as in lemma~\ref{lem-prune}, and it is a homotopy equivalence by that same lemma. \qed

\begin{defn}
\label{defn-tensormix}
We construct a map
\begin{eqnarray*}
(C\boxtimes D)\otimes(E\boxtimes F) & \lra & (C\otimes E)\boxtimes(D\otimes F)
\end{eqnarray*}
by composing the following:
\[
\xymatrix{
{\rule[-8mm]{0mm}{8mm}(C\boxtimes D)\otimes(E\boxtimes F)} \ar[d] \\
{\rule{0mm}{10mm}\holimsub{\twosub{\twosub{W=U\times V}{P_1=J_1\times K_1}}{P_2=J_2\times K_2}}\!\!\!\!
(C\otimes D)(W,W\smin P_1) \otimes (E\otimes F)(W,W\smin P_2)} \ar[d] \\
{\rule{0mm}{6mm}(C\otimes E)\boxtimes(D\otimes F)}
}
\]
The second arrow is a right inverse for the chain map in proposition~\ref{prop-boxes}. The first
arrow is induced by isomorphisms
\[
\xymatrix{
\big(C(U,U\smin J_1)\otimes D(U,U\smin J_2)\big)\otimes \big(E(V,V\smin K_1)\otimes F(V,V\smin K_2)\big)
\ar[d]^\cong \\
\big((C\otimes E)(W,W\smin P_1)\big)\otimes\big((D\otimes F)(W,W\smin P_2)\big)
}
\]
where $W=U\times V$ and $P_i=J_i\times K_i$ for $i=1,2$. We describe the composite map informally
as $\varphi\otimes\psi\mapsto \varphi\bar\otimes\psi$, where $\varphi$ and $\psi$ are chains in
$C\boxtimes D$ and $E\boxtimes F$, respectively.
\end{defn}

\begin{defn}
\label{defn-productsym} Assume $C=D$ and $E=F$ in definition~\ref{defn-tensormix}. We construct a map
\begin{eqnarray*}
(C\boxtimes C)^{h\ZZ/2}\otimes(E\boxtimes E)^{h\ZZ/2} & \lra & ((C\otimes E)\boxtimes(C\otimes E))^{h\ZZ/2}
\end{eqnarray*}
by composing the following:
\[
\xymatrix{
{\rule[-8mm]{0mm}{8mm}(C\boxtimes C)^{h\ZZ/2}\otimes(E\boxtimes E)^{h\ZZ/2}} \ar[d] \\
{\rule{0mm}{10mm}\bigg(\holimsub{\twosub{\twosub{W=U\times V}{P_1=J_1\times K_1}}{P_2=J_2\times K_2}}\!\!\!\!
(C\otimes E)(W,W\smin P_1) \otimes (C\otimes E)(W,W\smin P_2)\,\,\bigg)^{h\ZZ/2}} \ar[d] \\
{\rule{0mm}{6mm}\big((C\otimes E)\boxtimes(C\otimes E)\big)^{h\ZZ/2}}
}
\]
(details as in definition~\ref{defn-tensormix}, with superscripts $h\ZZ/2$ added on where appropriate).
We describe the composite map informally
as $\varphi\otimes\psi\mapsto \varphi\bar\otimes\psi$, where $\varphi$ and $\psi$ are chains in
$(C\boxtimes C)^{h\ZZ/2}$ and $(E\boxtimes E)^{h\ZZ/2}$, respectively.
\end{defn}

\begin{expl} In particular, suppose that
$(C,\varphi)$ and $(E,\psi)$ are ``symmetric objects'' (not necessarily
Poincar\ee) of dimensions $m$
and $n$ respectively, in $\sC_X$ and $\sC_Y$ respectively, so that $\varphi$ is
an $m$-cycle in $(C\boxtimes C)^{h\ZZ/2}$ and $\psi$ is an $n$-cycle in
$(D\boxtimes D)^{h\ZZ/2}$. Then we have
\[  (C\otimes D,\varphi\bar\otimes\psi), \]
a symmetric object of dimension $m+n$ in $\sC_{X\times Y}$.
\end{expl}

\begin{expl}
\label{expl-productformula} Let $X$ and $Y$ be compact ENRs. Let $C$ be the functor
taking an open $U\subset X$ to the singular chain complex of $U$. Let $E$ be the functor
taking an open $V\subset Y$ to the singular chain complex of $V$. Let $F$ be the functor
taking an open $W\subset X\times Y$ to the singular chain complex of $W$. By all the above, we
have the following diagram
\[
\xymatrix@C=80pt{
C(X)\otimes E(Y) \ar[d]^\simeq \ar[r]^{\nabla\otimes\nabla} &
(C\boxtimes C)^{h\ZZ/2}\otimes (E\boxtimes E)^{h\ZZ/2} \ar@{..>}[d] \\
F(X\times Y) \ar[r]^{\nabla} & (F\boxtimes F)^{h\ZZ/2}~.
}
\]
which is commutative up to a chain homotopy. (We could be more precise about that by
specifying a ``contractible choice'' of such chain homotopies.) The dotted arrow is given as in
definition~\ref{defn-productsym} by $\varphi\otimes\psi \mapsto \varphi\bar\otimes\psi$.
We leave it to the reader to establish the homotopy commutativity.
\end{expl}

\section{Duality and decomposability}
\label{sec-decomposable}
We turn to a discussion of duality, first in $\sC$, then in $\sD$ and then in $r\sD$.
Let $C$ and $D$ be objects of $\sC$ which admit duals $C^{(-*)}$ and
$D^{(-*)}$. (In other words, the functors $E\mapsto H_0(C\boxtimes E)$
and $E\mapsto H_0(D\boxtimes E)$ on $\sH\sD/\sH\sD''$ are co-representable.)
Fix nondegenerate $0$-cycles
\[ \varphi\in C\boxtimes C^{(-*)}~,\quad \psi\in D\boxtimes D^{(-*)}~. \]
Let $f\co C\to D$ be any morphism. We can assume that $D^{(-*)}$ and $C^{(-*)}$ are
in $\sC''$. Choose a morphism $g\co D^{(-*)}\to C^{(-*)}$
such that $g_*(\psi)\in D\boxtimes C^{(-*)}$ is homologous to
$f_*(\varphi)$. Then choose a $1$-chain $\zeta\in D\boxtimes C^{(-*)}$ such that
$d\zeta=g_*(\psi)-f_*(\varphi)$. Now for every $E$ in $\sC$ the square
\[
\xymatrix{
\hom(C^{(-*)},E) \ar[d]_{\textup{slant with }\varphi}
\ar[r]^{g_*} & \hom(D^{(-*)},E) \ar[d]^{\textup{slant with }\psi} \\
C\boxtimes E  \ar[r]^{f_*} & D\boxtimes E
}
\]
is homotopy commutative (slant with $\zeta$ provides a homotopy), and the vertical
arrows are homology equivalences. Writing $\ko(g)$ and $\ko(f)$ for the mapping cone
of $g$ and $f$ respectively, we have homotopy cofiber sequences
\[
\xymatrix{
\hom(\ko(g),E) \ar[r] & \hom(C^{(-*)},E) \ar[r]^{g_*} & \hom(D^{(-*)},E)~,
}
\]
\[
\xymatrix{
C\boxtimes E \ar[r]^{f_*} & D\boxtimes E \ar[r] & \ko(f)\boxtimes E~.
}
\]
Our homotopy commutative square therefore implies that the functor
$E\mapsto H_0(\ko(f)\boxtimes E)$ is again co-representable, with representing
object $\Sigma^{-1}\ko(g)$. In particular the identity class in
$H_0\hom(\ko(g),\ko(g))$ corresponds to some nondegenerate class
$[\lambda]\in H_1(\ko(f)\boxtimes\ko(g))$. This construction of $[\lambda]$
shows also that, for every open $U\in\sO(X)$ and closed $K\subset X$ contained in $U$,
we have a commutative diagram
\[
\xymatrix@C=70pt{
\vdots  & \vdots  \\
{H^{j+1}\ko(f)(U,U\smin K)} \ar[u] \ar[r]^{\textup{slant with }\lambda} & H_j\ko(g)(U) \ar[u] \\
{H^jC(U,U\smin K)} \ar[u] \ar[r]^{\textup{slant with }\varphi} & {H_jC^{(-*)}(U)} \ar[u]\\
{H^jD(U,U\smin K)} \ar[u]_{f^*} \ar[r]^{\textup{slant with }\psi} & {H_jD^{(-*)}(U)} \ar[u]^{g_*} \\
{H^j\ko(f)(U,U\smin K)} \ar[u] \ar[r]^{\textup{slant with }\lambda} & H_{j+1}\ko(g)(U) \ar[u] \\
\vdots \ar[u] & \vdots \ar[u]
}
\]
with exact columns. This leads us to the following conclusion.

\begin{lem}
\label{lem-dualexist}
Objects of $\sD$ have duals which are again in $\sD$.
For $E$ and $F$ in $\sD$, an element $[\lambda]\in H_0(E\boxtimes F)$ is
nondegenerate if and only if, for all open $U$ in $X$,
the slant product with $[\lambda]$ is an isomorphism
\[  \colim_K~H^jE(U,U\smin K) \la F(U) \]
(where $K$ runs through the closed subsets of $X$ which are contained in $U$).
\end{lem}

\proof By the preceding discussion, it is enough to show that an
object $E$ of $\sD$
constructed as in example~\ref{expl-guide} from a map
$f\co \Delta^k\to X$ has a dual $F$ in $\sC$, with nondegenerate $[\lambda]\in H_0(E\boxtimes F)$
say, that $F$ is again decomposable, and that the slant product with $\lambda$
is an isomorphism
\[  \colim_K~H^jE(U,U\smin K) \la F(U) \]
for all $U\in\sO(X)$. By example~\ref{expl-can2} and
proposition~\ref{prop-nondeg}, all that is true.
\qed

\begin{cor}
\label{cor-dualexist}
Objects of $r\sD$ have duals which are again in $r\sD$.
For $E$ and $F$ in $r\sD$, an element $[\lambda]\in H_0(E\boxtimes F)$ is
nondegenerate if and only if, for all open $U$ in $X$,
the slant product with $[\lambda]$ is an isomorphism
\[  \colim_K~H^jE(U,U\smin K) \la F(U)~. \]
\end{cor}

\proof Let $E$ be an object in $r\sD$. We can assume that $E$ is in $r\sD'$ and that
$E$ admits a ``complement'' $E'$~, also in $r\sD'$, so that $E''=E\oplus E'$ belongs to
$\sD$. Now $E''$ admits a dual, say $F''$ in $\sD'$, coupled to $E''$ by means of
\[  [\lambda'']\in H_0(E''\boxtimes F''). \]
The retraction map $q:E''\to E''$ (via $E$) has a dual $p:F''\to F''$, so that
$(q\otimes\id)_*[\lambda'']=(\id\otimes p)_*[\lambda'']\in H_0(E''\otimes F'')$.
As $q$ is idempotent, $p$ is idempotent up to homotopy. We can now produce a splitting
$F''\simeq F\oplus F'$. Namely, for open $U$ in $X$ we let $F(U)$ be the homotopy colimit of
\[
\xymatrix{F''(U) \ar[r]^p & F''(U) \ar[r]^p &  F''(U) \ar[r]^p & \cdots
}
\]
and we let $F'(U)$ be the homotopy colimit of
\[
\xymatrix@C=35pt{F''(U) \ar[r]^{\id-p} & F''(U) \ar[r]^{\id-p} &  F''(U) \ar[r]^{\id-p} & \cdots
.}
\]
Then it is clear that $F$ and $F'$ belong to $r\sD'$ and
\begin{eqnarray*} H_0(E''\boxtimes F'') & \cong & H_0(E\boxtimes F)\oplus
H_0(E'\boxtimes F')\oplus H_0(E'\boxtimes F)\oplus H_0(E\boxtimes F').
\end{eqnarray*}
The equation $(q\otimes\id)_*[\lambda'']=(\id\otimes p)_*[\lambda'']$
shows that $[\lambda'']$ lives in
\[  H_0(E\boxtimes F)\oplus
H_0(E'\boxtimes F')~. \]
Write $[\lambda'']=[\lambda]\oplus[\lambda']$ with $[\lambda]\in H_0(E\boxtimes F)$
and $[\lambda']\in H_0(E'\boxtimes F')$. It is straightforward to show that
$[\lambda]\in H_0(E\boxtimes F)$ is nondegenerate \emph{and} that the
slant product with $[\lambda]$ is an isomorphism
\[  \colim_K~H^jE(U,U\smin K) \la F(U)~, \]
because $[\lambda'']$ has the analogous properties. The universal property of $[\lambda]$ now
implies that, for arbitrary $F^\sharp$ in $r\sD$ with $[\lambda^\sharp]\in H_0(E\boxtimes F^\sharp)$,
the element $[\lambda^\sharp]$ is nondegenerate if and only if the slant product
with it is an isomorphism
\[  \colim_K~H^jE(U,U\smin K) \la F^\sharp(U)~. \]
\qed

\begin{lem} The rule $X\mapsto \sD_X$
is a covariant functor, preserving duality.
\end{lem}

\proof Let $f\co X\to Y$ be a map (between locally compact separable Hausdorff spaces).
It is clear that, for $C$ in $\sD_X$, we have $f_*D$ in $\sD_Y$.
For $C$ and $D$ in $\sD_X$, there is a specialization map
\[  C\boxtimes D \la f_*C\boxtimes f_*D. \]
We need to show that this takes nondegenerate classes in $H_0(C\boxtimes D)$ to
nondegenerate classes in $H_0(f_*C\boxtimes f_*D)$. In the case where $X$ and
$Y$ are both compact, this follows from the nondegeneracy criterion
given in lemma~\ref{lem-dualexist}. (For $U$ open in $Y$,
every compact subset of $f^{-1}(U)$ is contained in some $f^{-1}(K)$ for compact
$K\subset U$.) Finally, because of finiteness condition (ii) in definition~\ref{defn-sheaflist},
it is easy to reduce to a situation where $X$ and $Y$ are both compact.
\qed

\begin{cor} The rule $X\mapsto \sD_X$
is a covariant functor, preserving duality. \qed
\end{cor}

\begin{prop} The tensor product $\sD_X\times \sD_Y\la \sD_{X\times Y}$ is compatible with duality.
\end{prop}

\proof Let $C,E$ be objects of $\sD_X$ and let
$D,F$ be objects of $\sD_Y$. Let $[\lambda]\in H_0(C\boxtimes E)$ and $[\mu]\in H_0(D\boxtimes F)$.
Then we have $[\lambda\bar\otimes\mu]\in H_0((C\otimes D)\boxtimes(E\otimes F))$.
What we have to show is that if $[\lambda]$ and $[\mu]$ are nondegenerate, then $[\lambda\bar\otimes\mu]$
is nondegenerate. This is a statement about the three triangulated categories $\sH\sD'_X$, $\sH\sD'_Y$
and $\sH\sD'_{X\times Y}$. For each of the three triangulated categories we know that duality
preserves exact triangles. We also know that the $\boxtimes$ product preserves exact triangles when
one input variable is fixed. Hence, using repeated five lemma arguments, we can easily reduce the claim
about the nondegeneracy of $[\lambda\bar\otimes\mu]$ to the special case where $C$ and $D$ are among the
standard generators of $\sD_X$ and $\sD_Y$, respectively. That is, $C$ is weakly equivalent to the
object of $\sD_X$ obtained from some map $f\co \Delta^k\to X$ by the method of example~\ref{expl-guide}, and $D$
is weakly equivalent to the object of $\sD_X$ obtained from some map $g\co \Delta^\ell\to Y$ by the same method.
(It seems better to say ``weakly equivalent'' rather than ``equal'' because we might want to apply the
resolution procedure of lemma~\ref{lem-resolve} to obtain objects in $\sD'_X$ and $\sD'_Y$~, respectively.)
In that case, we also have a clear idea what $E$ and $F$ are, and what $[\lambda]$ and $[\mu]$ are. Namely,
$E$ is (up to desuspensions) the quotient of $C$ by its ``boundary'' (the object obtained from
$f|\partial\Delta^k$ by the method of example~\ref{expl-guide}), and $F$ is (up to desuspensions) the
quotient of $D$ by its ``boundary'' (the object obtained from $g|\partial\Delta^\ell$ by the method of
example~\ref{expl-guide}). Also, $[\lambda]$ can be described as the class of $\nabla(\omega_k)$ where
$\omega_k$ is a relative fundamental cycle for the manifold-with-boundary $\Delta^k$, and
$[\mu]$ can be described as the class of $\nabla(\omega_\ell)$ where
$\omega_\ell$ is a relative fundamental cycle for the manifold-with-boundary $\Delta^k$. \newline
Next, $C\otimes D$ can be identified with the object obtained from
\[ f\times g\co \Delta^k\times\Delta^\ell\lra X\times Y \]
by the method of example~\ref{expl-guide}. Also $E\otimes F$ can be identified (up to
desuspensions) with the quotient of $C\otimes D$ by its ``boundary'', which is the object obtained
from $f\times g$ restricted to $\partial(\Delta^k\times\Delta^\ell)$ by the method of example~\ref{expl-guide}.
Now example~\ref{expl-productformula} implies that $[\lambda\bar\otimes\mu]$ can be described as the class of
$\nabla(\omega_k\times\omega_\ell)$. Since $\omega_k\times\omega_\ell$ is a relative fundamental cycle
for $\Delta^k\times\Delta^\ell$, this implies (with examples~\ref{expl-can1},~\ref{expl-can2} and
proposition~\ref{prop-nondeg}) that $[\lambda\bar\otimes\mu]$ is indeed nondegenerate. \qed

\section{The excisive signature}

\begin{lem}
\label{lem-homotopyinvariant}
The functor $X\mapsto \bL^\bullet(\sD_X)$ is
homotopy invariant.
\end{lem}

\proof It is enough to show that the maps $X\to X\times[0,1]$ given
by $x\mapsto(x,0)$ and $x\mapsto(x,1)$ induce the same homomorphisms
\[ \pi_*\bL^\bullet(\sD_X)\lra \bL^\bullet(\sD_{X\times[0,1]}).
\]
That is easily done by using the 1-dimensional manifold with
boundary $[0,1]$ and the corresponding SAPC in $\sC_{[0,1]}$, and
tensor product with that, to produce appropriate bordisms. \qed

\begin{thm}
\label{thm-excisive}
The functor $X\mapsto \bL^\bullet(\sD_X)$ is excisive.
In detail:
\begin{itemize}
\item[(i)] For open $U,V$ subset $X$ with $U\cup V=X$, the commutative square
of inclusion-induced maps
\[ \CD
\bL^\bullet(\sD_{U\cap V}) @>>> \bL^\bullet(\sD_U) \\
@VVV   @VVV  \\
\bL^\bullet(\sD_{V}) @>>> \bL^\bullet(\sD_X).
\endCD
\]
is homotopy (co)cartesian~;
\item[(ii)] For a finite or countably infinite disjoint union $X=\coprod X_\alpha$~,
the inclusions $X_\alpha\to X$ induce a (weak) homotopy equivalence
\[  \bigvee_\alpha\bL^\bullet(\sD_{X_\alpha})
\lra \bL^\bullet(\sD_X)~. \]
\end{itemize}
\end{thm}

\proof Excision property (ii) is a straightforward consequence of finiteness
property (ii). With corollary~\ref{cor-transverse}, the proof of excision property (i)
can be given using a mechanism which is very nicely abstracted in a paper by Vogel
\cite[1.18, 6.1]{Vogel}. \qed

\bigskip
\begin{thm}
\label{thm-Rothenberg}
The relative homotopy groups of the inclusion $\bL^\bullet(\sD_X)\to \bL^\bullet(r\sD_X)$
are vector spaces over $\ZZ/2$.
\end{thm}

\proof Let $K_0(\sD)$ be the Grothendieck group of $\sD$ (with one generator $[C]$ for each object $C$, a
relation $[C]\sim[D]$ if $C$ and $D$ are weakly equivalent, and a relation $[C]-[D]+[E]=0$ for
every short exact sequence $C\to D\to E$ in $\sD$). Define the Grothendieck group of $r\sD$ similarly. Let
\[  \tilde K_0 \]
be the cokernel of the inclusion-induced map $K_0(\sD)\to K_0(r\sD)$. The group $\ZZ/2$ acts on
this by means of (degree 0) duality.
The long exact sequence of homotopy groups of the inclusion map $\bL^\bullet(\sD_X)\to \bL^\bullet(r\sD_X)$
can be described as a ``Rothenberg'' sequence:
\[ \cdots\to  L^n(\sD_X)\to L_n(r\sD_X) \to
\widehat H^n(\ZZ/2;\tilde K_0) \to L^{n-1}(\sD_X)\to L^{n-1}(r\sD_X)\to \cdots
\]
where $\widehat H^*$ denotes Tate cohomology. \qed

\begin{rem} Lemma~\ref{lem-homotopyinvariant} and theorems~\ref{thm-excisive} and~\ref{thm-Rothenberg}
have analogues for quadratic $L$-theory which can be proved in the same way.
\end{rem}

\medskip
If our locally compact Hausdorff separable space $X$ is an ENR, then
\[  \bL^\bullet(\sD_X) \simeq X_+\wedge \bL^\bullet(\sD_\pt)= X_+\wedge \bL^\bullet(\ZZ).\]
This follows from homotopy invariance and the two excision
properties by the standard arguments going back to Eilenberg and
Steenrod. Here $\bL^\bullet(\ZZ)$ is the symmetric $L$-theory
spectrum of the \emph{ring} $\ZZ$ (with the trivial involution). 

\begin{rem} If $X$ is the polyhedron of a simplicial complex
$\bL^\bullet(\sD_X)$ has the homotopy type of the spectrum 
$\bL^\bullet(\ZZ,X)$ constructed in \cite[\S10]{RanickiTop} from the 
$(\ZZ,X)$-category of \cite{RanickiWeiss1} endowed with a chain duality. 
See also \cite{Woolf}, \cite{RanickiWeiss2}, \cite{LauresMcclure}.
\end{rem}

If $X$ is a compact oriented topological $n$-manifold with boundary
$\partial X$, then the
identity map $X\to X$ determines by example~\ref{expl-guide} and example~\ref{expl-can2}
an $n$-dimensional SAP pair in $r\sD_X$ (with boundary in $r\sD_{\partial X}$)
which in turn determines an element in
\[
\pi_n(\bL^\bullet(r\sD_X),\bL^\bullet(r\sD_{\partial X})) \, \cong_{\ZZ[1/2]} \,
\pi_n(\bL^\bullet(\sD_X),\bL^\bullet(\sD_{\partial X}))\, \cong \, H_n(X,\partial X;\bL^\bullet(\ZZ)). \]

\begin{defn}
\label{defn-excisignature}
This element in $\pi_n(\bL^\bullet(r\sD_X),\bL^\bullet(r\sD_{\partial X}))$
is the excisive signature of $(X,\partial X)$.
\end{defn}

\begin{rem} If $X$ is triangulable, then we can regard the excisive signature of $(X,\partial X)$
as an element of $\pi_n(\bL^\bullet(\sD_X),\bL^\bullet(\sD_{\partial X}))$. In fact the
excisive signature of a compact topological manifold $X$ with boundary, not necessarily triangulable,
can \emph{always} be regarded as
an element of
\[ \pi_n(\bL^\bullet(\sD_X),\bL^\bullet(\sD_{\partial X}))~. \]
This follows easily from the fact that $X\times I^n$ for sufficiently large $n$ admits a handle decomposition.
Unfortunately the proof of that fact (existence of handle decomposition) given e.g. in \cite{KirbySiebenmann}
is hard and uses ideas which are quite closely related to Novikov's original proof of the topological
invariance of Pontryagin classes. For this reason we do not wish to use the ``handle decomposition'' argument.
We have already avoided it by introducing $r\sD_X$ and proving theorem~\ref{thm-Rothenberg}.
\end{rem}

\begin{rem} Let $f:Y\to X$ be a degree 1 normal map of closed $n$-dimensional topological manifolds. By
example~\ref{expl-guide} and example~\ref{expl-can2}, the map $f:Y\to X$ and the identity map
$\id:X\to X$ determine two $n$-dimensional SAP objects $(C(f),\varphi)$ and $(C(\id),\psi)$ in
$r\sD_X$ (and even in $\sD_X$, by the previous remark).
The map $f$ induces a chain map $C(f)\to C(\id)$ which respects the symmetric structures,
so that there is a splitting up to weak equivalence in $r\sD_X$ or $\sD_X$~,
\[  (C(f),\varphi)~\simeq~ (C(\id),\psi)\oplus(K,\zeta)~. \]
We expect that the nondegenerate symmetric structure $\zeta$ on $K$ has a canonical
refinement to a (nondegenerate) quadratic structure, determined by the bundle data which come with
the normal map $f$.
\end{rem}

\section{The Poincar\ee dual of the excisive signature}
There is a rational homotopy equivalence
\[ \bL^\bullet(\ZZ)\simeq_{\QQ} \bigvee_{i\ge 0} S^{4i}\wedge \bH\QQ~, \]
unique up to homotopy. For a compact oriented topological $n$-manifold $X$
with boundary, the Poincar\ee dual of the ``rationalized'' excisive signature
of $(X,\partial X)$ is therefore a class in
\[  \bigoplus_{i\ge 0} H^{4i}(X;\QQ) . \]
We shall show that it is a characteristic class
associated with the topological tangent bundle of $X$, a bundle
with structure group $\TOP(n)$.

\begin{lem} The suspension isomorphism
\[ H_0(\pt;\bL^\bullet(\ZZ))\la  H_1(I,\partial I;\bL^\bullet(\ZZ)) \]
takes the unit $1$ to the excisive signature of $(I,\partial I)$. \qed
\end{lem}

\begin{prop}
\label{prop-suspension}
Let $X$ be a compact oriented topological $n$-manifold $X$ with boundary and let $Y=X\times [0,1]$,
so that $Y/\partial Y\cong \Sigma(X/\partial X)$.
The suspension isomorphism
\[  H_n(X,\partial X;\bL^\bullet(\ZZ))\otimes \ZZ[1/2] \la H_{n+1}(Y,\partial Y;\bL^\bullet(\ZZ))\otimes\ZZ[1/2] \]
takes the excisive signature of $(X,\partial X)$ to the excisive signature of $(Y,\partial Y)$.
\end{prop}

\proof This follows from the previous lemma and the product formula in example~\ref{expl-productformula}.  \qed

\begin{prop}
\label{prop-collapse}
Let $X$ be a compact oriented topological $n$-manifold $X$ with boundary, $Y\subset X$
a compact codimension zero submanifold with locally flat boundary, $Y\cap \partial X=\emptyset$.
Then, under the homomorphism
\[ H_n(X,\partial X;\bL^\bullet(\ZZ))\otimes \ZZ[1/2] \la H_n(Y,\partial Y;\bL^\bullet(\ZZ))\otimes \ZZ[1/2] \]
induced by the quotient map $X/\partial X\to Y/\partial Y$, the excisive signature
of $(X,\partial X)$ maps to the excisive signature of $(Y,\partial Y)$.
\end{prop}

\proof This is a consequence of the naturality of $\nabla$ in example~\ref{expl-can1}. \qed

\begin{prop}
\label{prop-signaturedualized}
The Poincar\ee dual of the (rationalized) excisive signature of a compact oriented manifold
with boundary is a characteristic class $\Lambda$ for euclidean bundles, evaluated on the tangent (micro)bundle
of the manifold. (The characteristic class $\Lambda$ is defined for euclidean bundles on compact ENRs,
and is invariant under stabilisation, i.e., replacing a euclidean bundle $E\to Y$ by $E\times\RR\to Y$.)
\end{prop}

\begin{rem} While the construction of $\Lambda$ as such is elementary, we use a technical
fact from geometric topology to show that $\Lambda(TY)$ is Poincar\ee dual to $\sigma(Y,\partial Y)$
in the case where $Y$ is a compact $n$-manifold. This fact is the existence of stable normal bundles
for embeddings of topological manifolds \cite{Hirsch2}.
\end{rem}

\proof[Proof of proposition~\ref{prop-signaturedualized}]
 Let $Y$ be a finite simplicial complex and let $E\to Y$ be a bundle on $Y$ with fibers homeomorphic to $\RR^k$.
We would like to find a compact topological $n$-manifold $X$ for some $n$, and a homotopy equivalence
$f\co Y\to X$ such that $f^*TX$ is isomorphic to $E\times\RR^{n-k}\to Y$, a stabilized version of $E\to Y$.
Assuming that a sufficiently canonical choice of such an $X$ and $f$ can be made, we may then \emph{define} the
characteristic class associated with $E$ on $Y$ to be $f^*$ of the Poincar\ee dual of the excisive signature of
$(X,\partial X)$. \newline
For the first step of this program, we choose an embedding $Y\to \RR^\ell$ which is linear on each
simplex of $Y$. Let $Y_r$ be a regular neighborhood of $Y$ in $\RR^\ell$. Choose an extension
of $E$ to a euclidean bundle $E_r\to Y_r$. Let $X\to Y_r$ be the bundle of $(k+1)$-disks on $Y_r$
obtained from $E_r\to Y_r$ by fiberwise one-point-compactification, followed by fiberwise join with a point.
Then $X$ is a compact oriented manifold of dimension $\ell+k+1$. Let $f\co Y\to X$ be the composition
of the inclusion $Y\to Y_r$ with any section of $X\to Y_r$. It is clear that $f^*TX$ is
identified with $E\times\RR^{\ell+1}\to Y$. \newline
We now define, as promised,
\[  \Lambda(E\to Y) = f^*(u_X)\in \bigoplus_{i\ge 0} H^{4i}(Y;\QQ) \]
where $u_X\in \bigoplus_{i\ge 0} H^{4i}(X;\QQ)$ is the Poincar\ee dual of the rationalized excisive
signature $\sigma(X,\partial X)$, in other words
\[   u_X\cap [X,\partial X]= \sigma(X,\partial X)~. \]
From proposition~\ref{prop-suspension} we deduce that this is well defined, i.e., independent of the
choice of an $\ell$ and an embedding $Y\to \RR^\ell$. (More precisely proposition~\ref{prop-suspension}
gives us the permission to make $\ell$ as large as we like, and for large $\ell$ any two embeddings
$Y\to \RR^\ell$ are isotopic.) From proposition~\ref{prop-collapse}, we deduce that $\Lambda$
is a characteristic class. Namely, suppose that we have euclidean bundles $E\to Y$ and $E'\to Y'$ and
a simplicial map $g\co Y\to Y'$ such that $g^*E'\cong E$. Taking $\ell$ large,
we can choose embeddings $Y\to \RR^\ell$ and $Y'\to \RR^\ell$, with regular neighborhoods $Y_r$ and
$Y'_r$~, in such a way that there is a codimension zero embedding $g_r\co Y_r\to Y'_r$ making the
following diagram homotopy commutative:
\[
\CD
Y @>\textup{inclusion}>> Y_r \\
@VV g V    @VV g_r V \\
Y' @>\textup{inclusion}>> Y'_r~.
\endCD
\]
Now proposition~\ref{prop-collapse} can be applied to the embedding $g_r$ and gives the desired conclusion,
that $g^*\Lambda(E')=\Lambda(E)$ in $\bigoplus_{i\ge 0} H^{4i}(Y;\QQ)$. \newline
Finally we can mechanically extend the definition of $\Lambda$ to obtain a characteristic class
defined for euclidean bundles on compact ENRs. Indeed let $Y$ be a compact ENR; then $Y$ is a
retract of some finite simplicial complex $Y_1$. Hence any euclidean bundle on $Y$ extends to
one on $Y_1$. We can evaluate the characteristic class $\Lambda$ there, and pull back to the
cohomology of $Y$. To show that this is well defined, use the following: if $Y$ is a retract
of a finite simplicial complex $Y_1$~, and also a retract of a finite simplicial complex $Y_2$,
then the union of $Y_1$ and $Y_2$ along $Y$ is again an ENR. See \cite{Hu}. \newline
This is not the end of the proof, because we still have to show that
\begin{eqnarray*}
\Lambda(TY)\cap(Y,\partial Y) & = &
\sigma(Y,\partial Y) \,\,\in \bigoplus_{i\ge0} H_{n-4i}(Y,\partial Y;\QQ)
\end{eqnarray*}
holds in the case where $Y$ is a compact $n$-manifold with boundary. To establish this, we choose first of all
a locally flat embedding $Y\to \RR^\ell$ for some $\ell$. This can be done by the method of
\cite[Ch.1,Thm.3.4]{Hirsch1}. In view of this we write $Y\subset \RR^\ell$. Increasing $\ell$ if necessary,
we may also assume \cite{Hirsch2} that $Y$ has a normal microbundle in $\RR^{\ell}$, and by \cite{Kister}
we may also assume that it has a normal bundle $N\to Y$ in $\RR^\ell$. Choose a neighborhood $Y_r$
of $Y$ in $\RR^\ell$ such that $Y$ is a retract of $Y_r$ and $Y_r$ is a codimension zero PL submanifold
of $\RR^\ell$. According to our definition of $\Lambda$, we now have to extend the euclidean bundle
$TY\to Y$ to a euclidean bundle $E_r\to Y_r$ (which is easy). Then we should replace $E_r\to Y_r$ by
$E_r\times\RR\to Y_r$, which completes to a disk bundle $X\to Y_r$, etc.; we then have to find $\sigma(X,\partial X)$
and pass to Poincar\ee duals.--- Altogether we now have an
embedding $Y\to X$ by composing
\[
\CD
Y @>\textup{incl.}>> Y_r @>\textup{incl.}>> E_r\times\RR @>\textup{incl.}>> X
\endCD
\]
where the second arrow is any section of the bundle projection $E_r\times\RR\to Y_r$. To complete the
proof, it suffices to show that $Y$ has a trivial normal disk bundle $X'$ in $X$, and to apply
proposition~\ref{prop-collapse} to the inclusion $X'\to X$. Here we note that the existence of a trivial
normal bundle (with fibers $\cong\RR^{\ell+1}$) implies the existence of a trivial normal disk bundle.
But it is clear that $Y$ has a normal bundle in $X$, identified with $N\times_Y(TY\times\RR)$,
and this is clearly trivial since already $N\times_YTY$ is trivial.
\qed

\begin{prop} On vector bundles, the characteristic class $\Lambda$ agrees with Hirzebruch's total
$\sL$-class.
\end{prop}

\proof Let $Y$ be a closed oriented topological manifold of
dimension $4i$. By construction, the scalar product $\langle
\Lambda_{4i}(TY),[Y]\rangle$ is equal to the image of $\sigma(Y)$
under the specialization (alias assembly) map
\[ H_{4i}(Y;\bL^\bullet(\ZZ))\otimes\QQ \to \pi_{4i}\bL^\bullet(\ZZ)\otimes\QQ~.  \]
In other words it is equal to the signature of $Y$. This holds in particular when $Y$
is smooth. As this property characterizes the $\sL$-class on vector bundles, we have
$\Lambda=\sL$ on vector bundles.
\qed

\end{document}